# Circumspheres of sets of $n+1$ random points in the $d$-dimensional Euclidean unit ball ($1 \leq n \leq d$)

Gérard Le Caër

*Institut de Physique de Rennes, UMR UR1-CNRS 6251, Université de Rennes I, Campus de Beaulieu, Bâtiment 11A, F-35042 Rennes Cedex, France*

*Electronic mail*: gerard.le-caer@univ-rennes1.fr

*Abstract*

In the $d$-dimensional Euclidean space, any set of $n+1$ independent random points, uniformly distributed in the interior of a unit ball of center $O$, determines almost surely a circumsphere of center $C$ and of radius $\Omega$ $(1 \leq n \leq d)$ and a $n$-flat $(1 \leq n \leq d-1)$. The orthogonal projection of $O$ onto this flat is called $O'$ while $\Delta$ designates the distance $O'C$. The classical problem of the distance between two random points in a unit ball corresponds to $n=1$. The focus is set on the family of circumspheres which are contained in this unit ball. For any $d \geq 2$ and any $1 \leq n \leq d-1$, the joint probability density function of the distance $\Delta \equiv O'C$ and of the circumradius $\Omega$ has a simple closed-form expression. The marginal probability density functions of $\Delta$ and of $\Omega$ are both products of powers and of a Gauss hypergeometric function. Stochastic representations of the latter random variables are described in terms of geometric means of two independent beta random variables. For $n = d \geq 1$, $\Delta$ and $\Omega$ have a Dirichlet distribution with parameters $(d, d^2, 1)$ while $\Delta$ and $\Omega$ are beta distributed. Results of Monte-Carlo simulations are in very good agreement with their calculated counterparts. The tail behavior of the circumradius probability density function has been studied by Monte-Carlo simulations for $2 \leq n = d \leq 9$, where all circumspheres are this time considered, regardless of whether or not they are entirely contained in the unit ball.

**I. INTRODUCTION**

Sets of $n+1$ points $A_i (i=1,..,n+1)$ are independently and uniformly distributed in the interior of a unit ball of center $O$ in the $d$-dimensional Euclidean space $\mathbb{R}^d$ for $d \geq 1$ and $1 \leq n \leq d$. The latter ball is denoted as $B_d$, where

$$B_d := \left\{ \mathbf{x} \in \mathbb{R}^d : \sum_{i=1}^{d} x_i^2 = \|\mathbf{x}\|^2 \leq 1 \right\}$$

and $\|\mathbf{x}\|$ is the Euclidean norm of $\mathbf{x}$. Each set of $n+1$ points determines almost surely a unique $n$-flat of $\mathbb{R}^d$, called $L_n$, for $d \geq 2$ and $1 \leq n \leq d-1$, an affine linear subspace of dimension $n$ which does not necessarily contain $O$ whose orthogonal projection onto it is $O'$ (figure 1). These $n+1$ points are almost never contained in any $(n-1)$-flat. For $n = d$, $O'$ coincides with $O$ (figure 1c for $n = d = 2$). We define then the distance $H \equiv OO'$ for $1 \leq n \leq d-1$. Further, each set of points $(A_1,..,A_{n+1})$ determines almost surely a unique circumsphere passing through them and their convex hull is almost surely a $n$-simplex (see for instance chapter 8 of [55]). Besides circumspheres, the specific hyperspheres we shall consider are either $(n-1)$ dimensional spheres of radius $r$, denoted as $S_r^{(n-1)}$, and $S_1^{(d-1)}$, the unit hypersphere of $\mathbb{R}^d$ of center $O$. The intersection of $L_n$ with $S_1^{(d-1)}$ is $S_r^{(n-1)}$ of center $O'$ (figure 1d for $d = 3, n = 2$). In the context of this work, circumspheres may be classified into three families, designated from now on respectively as $C_d^{(n)}$, $D_d^{(n)}$ and $E_d^{(n)}$ $(d \geq 1, 1 \leq n \leq d)$. Circumspheres of $C_d^{(n)}$ are entirely contained in $B_d$. Those of $D_d^{(n)}$ and $E_d^{(n)}$ cut $S_1^{(d-1)}$. Their centers are respectively inside and outside the unit $d$-ball. The radii range between zero and one, zero and two and zero and infinity for circumspheres of $C_d^{(n)}$, $D_d^{(n)}$ and $E_d^{(n)}$ respectively. The probabilities that a circumsphere belongs to $C_d^{(n)}$, $D_d^{(n)}$ or $E_d^{(n)}$ are denoted respectively as $P_d^{(n)}$, $Q_d^{(n)}$ and $R_d^{(n)}$ $\left( P_d^{(n)} + Q_d^{(n)} + R_d^{(n)} = 1 \right)$ (see figure 7a for $2 \leq n = d \leq 9$).

In the present article, we focus on circumspheres of the $C_d^{(n)}$ family. The probability $P_d^{(n)}$ (eq. 80, $d \geq 2$) was calculated by Affentranger [3]. To fix ideas, we note that $P_d^{(n)}$ is far smaller than 1 for $n > \sim 5$. It is larger than ~0.05 for $n < \sim 4$ and any $d \geq 5$. When $d$ increases from $n$ to infinity, $P_d^{(n)}$ increases for instance from 0.4 to ~0.6 for $n = 2$ (circles), from ~0.12 to ~0.3 for $n = 3$ (spheres) and from ~0.03 to ~0.13 for $n = 4$ (hyperspheres). A circumsphere of $C_d^{(n)}$ whose radius is $\omega$ is denoted henceforth either as $c_d^{(n)}(\omega)$ or simply as $c_d^{(n)}$ when radius specification is not needed. An example, shown in figure 1a, summarizes, for $1 \leq n \leq d-1$ and for circumspheres $c_d^{(n)}(\omega)$, the random variables (rv's) we are interested in. The plane of the figure is determined by the center $O$ of $B_d$, its orthogonal projection $O'$ onto $L_n$ and the circumcenter $C$ which lies in $L_n$ [7]. The line determined by $O'$ and $C$ cuts the circumsphere in $A$ and $B$ where $A$ is chosen to be closer to $O'$ than $B$ is and cuts $S_1^{(d-1)}$ in two points $D$ and $D'$ symmetric with respect to $O'$ ($D'$ is not shown). Figures 1c and 1d show schematic representations of the case of three points $A_1$, $A_2$ and $A_3$ $(n = 2)$ respectively in a 2D unit disc and in a 3D unit ball. The random variables considered in the present work are the lengths $\Delta \equiv O'C$ and $\Delta_C \equiv OC$ as well as the circumradius $\Omega \equiv CA \equiv CB$. In addition, we define two random variables $\Sigma = \Delta + \Omega$ and $T = 1 - \Sigma$. The rv $\Sigma$ has a simple geometrical meaning which appears in figure 1a. A non-degenerate triangle $O'CM$ $(\delta > 0, \omega > 0)$, where $M$ is any point of the circumsphere except $B$, yields $O'M < O'C + CM = \sigma$. If $M$ coincides with $B$, then the triangle becomes a segment and $O'B = \sigma$. Therefore, a circumsphere is contained in the unit ball if and only if $\sigma < r$. We notice, en passant, that the rv's $\Delta$, $\Delta_C$, $\Omega$ are also relevant for circumspheres which cut $S_1^{(d-1)}$. We shall derive, for any $d \geq 1$ and any $1 \leq n \leq d$ and for circumspheres of $C_d^{(n)}$, the joint probability density function (pdf) of the length $\Delta$ and of the circumradius $\Omega$ as well as the associated marginal pdf's.

For reasons discussed below, the case of two random points $(n=1)$, $A_1$ and $A_2$, uniformly distributed in the interior of $B_d$, is of particular interest (figure 1b). The circumsphere is a 0-sphere which consists of the pair $A_1, A_2$. It is entirely contained in $B_d$. The circumcenter $C$ is the middle of $A_1 A_2$, and the circumradius $\Omega$ is half of the length $A_1 A_2$. The 1-flat, $L_1$, determined by $A_1$ and $A_2$ is a line. Although it sounds somewhat artificial to consider the distance $A_1 A_2$ as a circumdiameter, this point of view allows us to use a unified approach to derive in section V the pdf's mentioned above for any $d \geq 1$ and any $1 \leq n \leq d$. This approach is not the most direct for $n=1$ but it differs from published ones. The pdf of the length $A_1 A_2$ has been calculated during the last hundred years by a variety of methods which yield different formal expressions [6, 15-16, 18-20, 24-25, 30, 36, 40, 44, 46, 48-49, 52, 58, 60-61, 63].

The pdf's of the distance between two points randomly distributed in the inside of spheres or of ellipsoids find numerous applications. For instance, they were shown to simplify the calculation of self-energies of spherically symmetric matter distributions interacting by means of radially symmetric two-body potentials [48-49]. These calculations were extended to ellipsoids as a first step towards convex bodies whose shapes deviate from spherical. García-Pelayo [24] derived the distance pdf in ellipsoids as an integral he applied to a study of the shape of the earth. Other physical applications of distance pdf's include in particular the use of double electron-electron resonance to study spherical aggregates with shell structure [34] and the field of wireless networks whose properties are strongly influenced by distances between nodes [45-46, 58]. Finally, it is worth mentioning the connection between distance pdf's and pdf's of random chord length in convex bodies, which depends on the considered secant randomness (see for instance [13]). The chord length pdf's apply for instance in the fields of neutronics and of reactor physics [19, 26, 35, 37, 41-42, 52, 57, 63]. Extensions of the previous problem include for instance the determination of the mean distance between a reference point and its $n$ th neighbour among a collection of $N$ points uniformly distributed in a hypersphere or in a hypercube of unit volumes in a $d$-dimensional Euclidean space [8]. Applying the results of the latter authors, Kowalski [38] gave a geometrical interpretation to a generalization of a distribution used to represent pion distribution in hadronic production models. Circumcircles and circumspheres play an important role in computational geometry. Domains of all kinds are meshed with Delaunay triangulations and Voronoi tessellations are constructed from them. Direct applications of circumspheres are much less common than the previous ones. An example is the analysis of protein-induced distortions in [4Fe-4S] atom clusters [23].

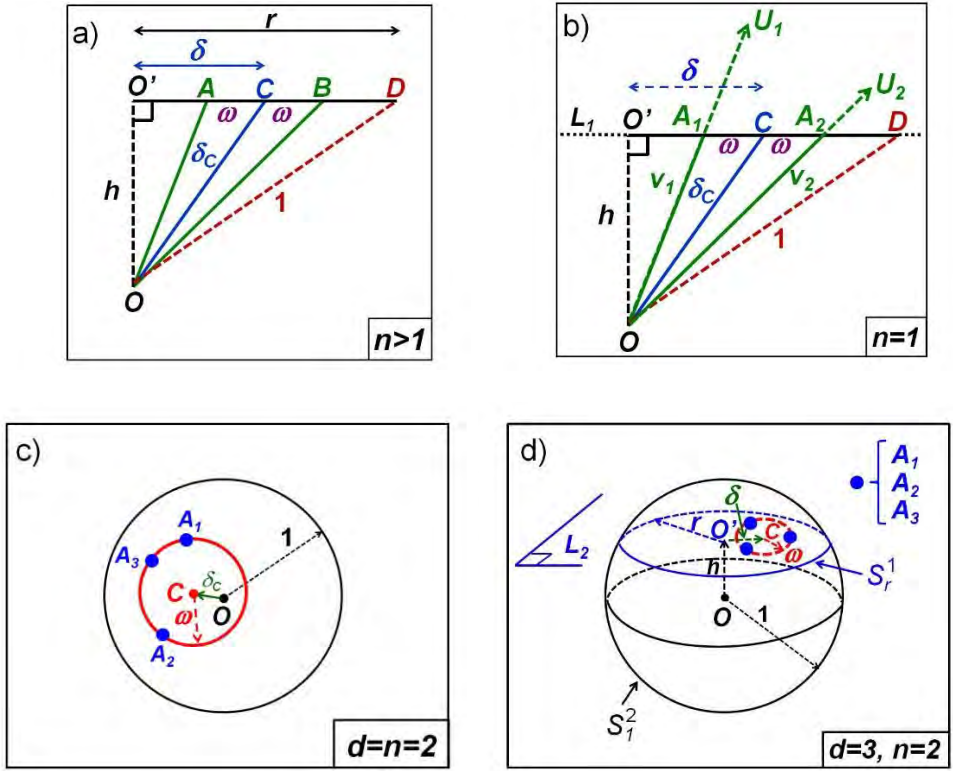

Figure 1: $O$ is the origin of the unit ball $B_d$, $O'$ its orthogonal projection onto the $n$-flat, $L_n$, determined by $n+1$ i.i.d. points $A_k$ $(k=1,..,n+1)$ uniformly distributed over $B_d$ $(1 \leq n \leq d-1)$. The circumsphere defined by this set of points has a center $C$ $(\delta = O'C, \delta_C = OC)$ and a radius $\omega$. The intersection of $L_n$ with the unit hypersphere of center $O$ is a hypersphere $S_r^{(n-1)}$ of center $O'$ and of radius $r$ cut by $O'C$ in $D$.

a) for $n>1$: $O'D$ cuts the circumsphere in $A$ and in $B$ and $\sigma = O'B = \delta + \omega$

b) For $n=1$, the two points $A$ and $B$ coincide respectively with $A_1$ and $A_2$ which determine a line $L_1$. These points are obtained from two i.i.d. unit vectors $U_1$ and $U_2$ and two i.i.d. rv's $V_1$ and $V_2$ (end of section IIIB).

An overall view is shown for the case of three points $A_1$, $A_2$ and $A_3$

c) in a unit disc ($O'$ coincides with $O$)

d) in a 3D unit ball where they determine a plane $L_2$ (for clarity the segment $OO'$ is chosen to be "vertical" without loss of generality).

The method we shall use to derive the joint distribution of $\Delta$ and $\Omega$ for circumspheres of $C_d^{(n)}$ is based:

1) on affine equivalence: Kingman [37] considered $n+1$ points $A_i (i=1,..,n+1)$ chosen at random within a convex body $K \subset R^d$ and the associated $n$-flat, $L_n$. As $K$ is here spherical, the sets $K \cap L_n$ are affinely equivalent for $n \geq 1$ and for different $L_n$ provided that the "volume" of $K \cap L_n$ is non-zero.

Thus, the random geometrical characteristics we investigate are first rescaled as follows. The intersection of $S_1^{(d-1)}$ with a given $L_n$ is $S_r^{(n-1)}$ whose radius $R = r$ (figures 1a and 1d) is chosen to become equal to 1 once rescaled. We define thus rescaled random variables $\Delta_r = \Delta/r$, $\Omega_r = \Omega/r$ and $T_r$ is defined as $T_r = (1-\Delta-\Omega)/r$ (sections III and IV).

2) on the results of Affentranger [2-3] which imply that the triplet $(\Delta_r, \Omega_r, T_r)$ has a Dirichlet distribution (section IVA).

3) on the known distribution $p_d^{(n)}(h)$ of the distance $H \equiv OO'$ [44] which yields the distribution $p_d^{(n)}(r)$ of the previous scaling factor $R = \sqrt{1-H^2}$ (figures 1a and 1c).

4) on an integration over $r$, with weights given by $p_d^{(n)}(r)$, of the conditional distribution of $\Delta$ and $\Omega$ given $R = r$ (sections III and IV).

First, we shall derive the joint distribution of the length $\Delta \equiv O'C$ and of the circumradius $\Omega$ for any $d \geq 1$ and any $1 \leq n \leq d$, as well as the mathematical form of the joint distribution of the length $\Delta_C \equiv OC$ and of $\Omega$. Second, we shall obtain the marginal distributions of $\Delta$, $\Delta_C$ and $\Omega$, their moments and simple stochastic representations of $\Delta$ and $\Omega$. Third, we shall calculate the probability $p\left(O' \notin c_d^{(n)}\right)$ that $O'$ is not contained in the interior of a circumsphere $c_d^{(n)}$. Fourth, we shall briefly discuss the asymptotic behaviour of the latter variables when $d \to \infty$ for a fixed $n$. Last, we shall report on results of a Monte-Carlo study of the tail behavior of the circumradius pdf, where all circumspheres are this time considered, irrespective of the fact whether or not they are entirely contained in the unit ball.

.

## II. SOME DEFINITIONS AND NOTATIONS

A $n$-simplex is the convex hull of $n+1$ affinely independent points. The vertices of the standard or unit $n$-simplex, $\Delta_n$, embedded in a hyperplane of $\mathbb{R}^{n+1}$, are the $n+1$ points:

$$\Delta_n : \left\{ v_1 = \overbrace{(1,0,...,0)}^{n+1}, v_2 = (0,1,...,0),..., v_{n+1} = (0,0,...,1) \right\} \quad (1)$$

A $n$-simplex is regular if the intervertex distances are all equal, say to $a_0$, with $a_0 = \sqrt{2}$ for the unit simplex. The circumradius of a regular $n$-simplex is obtained from a general relation satisfied by it [9]. Let a general point be at distance $a_k$ from vertex $k$, $(k=1,..,n+1)$, then:

$$(n+1)\left(\sum_{k=0}^{n+1} a_k^4\right) = \left(\sum_{k=0}^{n+1} a_k^2\right)^2 \qquad (2)$$

The circumradius of this regular $n$-simplex is thus $\left(a_0/\sqrt{2}\right) \times \sqrt{n/(n+1)}$ ([33] p. 257) which reduces to $\sqrt{n/(n+1)}$ for a unit $n$-simplex.

The beta function $B(q_1, q_2)$ is related to the Euler gamma function $\Gamma(q)$ by $B(q_1, q_2) = \Gamma(q_1)\Gamma(q_2)/\Gamma(q_1+q_2)$, where the parameters $q_1, q_2$ are here real and positive. The Pochhammer symbol $(a)_s$, where $s$ is a non-negative integer, is expressed as $(a)_s = \Gamma(a+s)/\Gamma(a) = a(a+1)..(a+s-1)$, $(a)_0 = 1$. The regularized incomplete beta function $I_x(a,b)$ is defined from the incomplete beta function by:

$$I_x(a,b) = \frac{\int_0^x t^{a-1}(1-t)^{b-1} dt}{B(a,b)} = \frac{x^a}{aB(a,b)} {}_2F_1(a, 1-b; a+1; x) \qquad (3)$$

(section 6.6 of [1]). It is formulated in terms of the Gauss hypergeometric function in the right-hand side of eq. 3.

Upper-case letters are used to denote random variables and lower-case letters for the values they take. The notation $A \triangleq B$ means that random variables $A$ and $B$ have the same distribution. Independent and identically distributed rv's are abbreviated as i.i.d..The joint probability density function of the continuous random variables $\Delta, \Omega$ and the pdf's of $X$, which represents any rv among the following seven $(\Delta, \Delta_C, \Omega, \Sigma, T, H, R)$ defined in the previous section, will be denoted respectively as $p_d^{(n)}(\delta, \omega)$ and $p_d^{(n)}(x)$. The use of a single notation to designate bivariate or univariate pdf's, whose functional forms differ, is not ambiguous as it is the nature of the random variables themselves which makes it clear which distributions $p_d^{(n)}(\ )$ are being dealt with. The conditional pdf of $X$ given $Y$ is denoted as usual $p_d^{(n)}(x|Y=y)$.

Besides the Gaussian distribution used in Monte-Carlo simulations (appendix A), the classical distributions considered in the present work are, the gamma distribution of shape parameter $q > 0$ and scale parameter $\theta > 0$, $X \sim \gamma(q, \theta)$, the beta distribution of shape parameters $q_1 > 0$ and $q_2 > 0$, $X \sim Be(q_1, q_2)$ and the Dirichlet distribution, $\boldsymbol{L}_k \sim Dir(\boldsymbol{q}_k)$.

Here $X$ and $L_k$ are respectively a random variable and a $k$-dimensional random vector while $\sim$ stands for "is distributed as". The Dirichlet distribution depends on a $k$-dimensional vector of positive parameters $q_k$. Further, it is defined on the unit $k-1$ simplex (eq. 1). The previous distributions, which are linked together, are discussed further in appendix B. The distribution of the product of two independent beta random variables is discussed in appendix C. It is used to derive stochastic representations of $\Delta$ and of $\Omega$ (appendices C and D).

Probability density functions and probabilities (table 1 for $n=2$, appendix A) were estimated for circumspheres of the $C_d^{(n)}$ family for various values of $n$ and of space dimensions $d$ from results of Monte-Carlo simulations, most often of $2.10^8$ circumspheres, with a method described in appendix A. Extreme value pdf's, which are determined by circumspheres of the $E_d^{(n)}$ family, were estimated too for $2 \leq n = d \leq 9$ (section VIII). The value $n=2$ (circumcircles) was selected for plots of some estimated pdf's shown in figures 2 to 6 and in figure 8. In these figures, points are placed at the midpoints of bins of size 0.001 and solid lines are drawn through them. The differences between simulated and calculated results are of the order of line thicknesses. Dotted vertical lines in figures 3b, 4, 5 and 6b are the asymptotic limits of the means of the variables (section VII) whose pdf's are shown in them.

### III. THE CASE OF TWO POINTS: $n=1$

The circumsphere $c_d^{(1)}$ of two points $A_1$ and $A_2$, which are independently and uniformly distributed in the interior of $B_d$, is a 0-sphere which consists of the pair $A_1$, $A_2$. All relevant characteristics are indicated in figure 1b. The simplest case, $d=1, n=1$, is first discussed.

### A. $d=1$

For $d=1$, two i.i.d. random variables $U_1$ and $U_2$, uniformly distributed over (-1,+1) yield two points, $A_1$ and $A_2$ respectively. The latter rv's, whose joint pdf is $p(u_1, u_2) = 1/4$, are then combined to give $Y_1 = (U_1 + U_2)/2$ and $Y_2 = (U_1 - U_2)/2$ with $p(y_1, y_2) = 1/2$, $y_1 + y_2 \in (-1, 1)$ and $y_1 - y_2 \in (-1, 1)$. The distance $\Delta$ and the circumradius $\Omega$ are respectively $|Y_1|$ and $|Y_2|$, where the condition $0 < \delta + \omega < 1$ is obviously obeyed. We define thus $T = 1 - \Delta - \Omega$. The joint pdf $p_1^{(1)}(\delta, \omega)$ is therefore:

$$p_1^{(1)}(\delta, \omega) = 2 \qquad \tau + \delta + \omega = 1 \qquad \delta, \omega, \tau \in (0,1) \tag{4}$$

The latter distribution is uniform over the standard 2-simplex. It is a Dirichlet distribution (appendix B) $Dir(1,1,1)$ of the triplet $(\Delta, \Omega, T)$ in agreement with eq. 20. As the distributions of $U_1$ and $U_2$ are symmetric, $U_k \triangleq -U_k$ $(k=1,2)$, then $Y_1 \triangleq Y_2$ and $\Delta \triangleq \Omega$. The common marginal distribution of $\Delta$ and $\Omega$ is derived either by a simple direct calculation or from the Dirichlet distribution $p_1^{(1)}(\delta, \omega)$ (eq. 89):

$$\begin{cases} p_1^{(1)}(\delta) = 2(1-\delta) \\ p_1^{(1)}(\omega) = 2(1-\omega) \end{cases} \qquad \delta, \omega \in (0,1) \tag{5}$$

in agreement with pdf's given in the next subsection by eq. 11 for $d=1$ as $_2F_1(1,1/2;2;1-\omega^2) = 2/(1+\omega)$ (eq. 7.3.1.125 of [50]).

**B.** $d \geq 2$

We consider more particularly line segments $B_d \cap L_1$ of length $2R$. The line segments $B_d \cap L_1$ are affinely equivalent for different $L_1$ provided that $R > 0$ [37]. Given the line $L_1$ and thus $R = r$, the coordinates of the points $A_1$ and $A_2$ on this secant, respectively $X_1$ and $X_2$, are then scaled in such a way that they range from 0 to 1. Kingman [37] (see also [54], p. 201, eq. 12.23) derived then the bivariate pdf of $X_1$ and $X_2$, $p(x_1, x_2)$ (eq. 6), which results from a Blaschke-Petkantschin type formula (section 7.2 of [55], see too section IVA below) applied to the case $n=1$. The volume of the simplex, whose vertices are just $A_1$ and $A_2$, reduces simply to $|x_1 - x_2|$. It is raised to the power $d-1$ in the previous formula, that is:

$$p(x_1, x_2) = \frac{d(d+1)}{2} \times |x_1 - x_2|^{d-1} \qquad x_1, x_2 \in (0,1) \tag{6}$$

We define two new random variables $U = X_1 + X_2 - 1$ and $Z = X_2 - X_1$ whose joint pdf is:

$$q(u,z) = \begin{cases} \dfrac{d(d+1)}{4} \times z^{d-1} & z \in (0, \min(1-u, 1+u)) \\ \dfrac{d(d+1)}{4} \times (-z)^{d-1} & z \in (\max(-1-u, -1+u), 0) \end{cases} \qquad u, z \in (-1,1) \tag{7}$$

The joint pdf of $\Delta_r = \Delta/r = |U|$ and of $\Omega_r = \Omega/r = |Z|$ is then $4q(\delta_r, \omega_r)$ (first line of eq. 7) with $\delta_r \in (0,1)$ and $\omega_r \in (0, 1-\delta_r)$. As $\delta_r + \omega_r \in (0,1)$, we define in addition, $T_r = 1 - \Delta_r - \Omega_r$ and we write finally:

$$f_d^{(1)}(\delta_r, \omega_r) = d(d+1)\omega_r^{d-1} \qquad \delta_r + \omega_r + \tau_r = 1 \qquad \delta_r, \omega_r, \tau_r \in (0,1) \tag{8}$$

The latter pdf is that of a Dirichlet distribution $Dir(1,d,1)$ of the triplet $(\Delta_r, \Omega_r, T_r)$ in agreement with eq. 20. As the pdf of $R$ will be shown to be $w_d^{(1)}(r) \propto r^{(d+2)}(1-r^2)^{(d-3)/2}$ (eq. 22, section IVB), the bivariate distribution $p_d^{(1)}(\delta, \omega)$ is obtained from $f_d^{(1)}(\delta_r, \omega_r)$:

$$p_d^{(1)}(\delta, \omega) \propto \omega^{d-1} \int_{\delta+\omega}^{1} w_d^{(1)}(r) \frac{dr}{r^{d+1}} \propto \omega^{d-1} \int_{\delta+\omega}^{1} 2r(1-r^2)^{\frac{d-3}{2}} dr \tag{9}$$

Finally:

$$p_d^{(1)}(\delta, \omega) = \frac{2^d}{\sqrt{\pi}} \times \frac{d^2 \Gamma(d/2)}{\Gamma((d+1)/2)} \times \omega^{d-1}\left(1-(\delta+\omega)^2\right)^{\frac{d-1}{2}} \qquad (\delta, \omega, \delta+\omega \in (0,1)) \tag{10}$$

The derivation of the marginal distribution $p_d^{(1)}(\omega)$ from $p_d^{(1)}(\delta, \omega)$ will not be reproduced here as the calculation of $p_d^{(n)}(\omega)$ from $p_d^{(n)}(\delta, \omega)$ is detailed in section VC for any $d \geq 1$ and is valid for any $1 \leq n \leq d$. The pdf $p_d^{(1)}(\omega)$ obtained from eqs 40 and 41 for $n=1$ (middle member of eq. 11) can then be rewritten in terms of the regularized incomplete beta function (eq. 3):

$$p_d^{(1)}(\omega) = \frac{2^d d^2 \Gamma(d/2)}{(d+1)\Gamma((d+1)/2)\sqrt{\pi}} \times \omega^{d-1} \times \left[(1-\omega^2)^{(d+1)/2} \,{}_2F_1\left(\frac{d+1}{2}, \frac{1}{2}; \frac{d+3}{2}; 1-\omega^2\right)\right] =$$

$$= 2^d d\omega^{d-1} I_{1-\omega^2}\left(\frac{d+1}{2}, \frac{1}{2}\right) \tag{11}$$

The right-hand side of eq. 11 is, as expected, identical with the density of the half length $\Omega = A_1 A_2/2$ given for instance by Hammersley [30] and Lord [40].

The cumulative distribution function $F_d(w)$ of the distance $W$ between a point uniformly distributed over the unit ball $B_d$ and its center $O$ is $F_d(w) = w^d$. The generation of $A_1$ and $A_2$ requires two i.i.d. random unit vectors, $U_1$ and $U_2$, which are uniformly distributed over the surface of the unit sphere in $\mathbb{R}^d$. Then two i.i.d. random variables $Z_1$ and $Z_2$, uniformly distributed over $(0,1)$, are used to write finally: $\boldsymbol{OA_1} \triangleq V_1 \boldsymbol{U_1}$ and $\boldsymbol{OA_2} \triangleq V_2 \boldsymbol{U_2}$ with $V_1 \triangleq (Z_1)^{1/d}$ and $V_2 \triangleq (Z_2)^{1/d}$ (figure 1b and appendix A). The distributions of the vectors $\boldsymbol{OC} = (V_1 \boldsymbol{U_1} + V_2 \boldsymbol{U_2})/2$ and $\boldsymbol{A_1 C} = (-V_1 \boldsymbol{U_1} + V_2 \boldsymbol{U_2})/2$ are identical because $\boldsymbol{U_k} \triangleq -\boldsymbol{U_k}$ ($k=1,2$) and because the rv's $(V_1, V_2, \boldsymbol{U_1}, \boldsymbol{U_2})$ are mutually independent. Thus, spherical symmetry imposes that the marginal distribution of the length $\Delta_C \equiv OC$ is identical with that of $\Omega$ for $n=1$, $\Delta_C \triangleq \Omega$ (eq. 11).

We note finally that the probability that the orthogonal projection $O'$ of $O$ onto $L_1$ does not belong to $A_1 A_2$ is simply obtained from $p(x_1, x_2)$ (eq. 6):

$$p\left(O' \notin c_d^{(1)}\right) = p\left(0 < X_1, X_2 < \frac{1}{2}\right) + p\left(\frac{1}{2} < X_1, X_2 < 1\right) = \frac{1}{2^d} \tag{12}$$

## IV. PRELIMINARIES FOR THE GENERAL CASE

### A. Joint distribution of $\Delta_r$ and of $\Omega_r$

We choose the case of three points $A_1, A_2, A_3$ in $\mathbb{R}^d$ ($n=2$) as detailed in [2] to exemplify the type of calculation performed in the general case by Affentranger [3]. The three points determine a unique 2-flat $L_2$. The volume element of $\mathbb{R}^d$ at a point $A_k$ is denoted as $dA_k$. By exterior multiplication, kept implicit in the notations, the Blaschke-Petkantschin formula writes ([2], p. 201 of [54]):

$$dA_1 dA_2 dA_3 = (2T)^{d-2} dA_1' dA_2' dA_3' dL_2 \tag{13}$$

Schneider and Weil (p. 272 of [51]) describe the common features of 'Blaschke–Petkantschin type' transformations which enables us to explain the structure of eq. 13. In their words, "the starting point of a transformation of 'Blaschke–Petkantschin type' is an integration over a product (possibly with one factor only) of measure spaces of geometric objects (points or flats as a rule), mostly homogeneous spaces with their invariant measures. Almost everywhere, the integration variable, which is a tuple of geometric objects, determines a new geometric object (for example, by span or intersection). We call this new object the 'pivot'. The initial integration is then decomposed into an outer and an inner integration. The outer integration space is the space of all possible pivots, with a natural measure; often it is a homogeneous space. For a given pivot, the inner integration space consists of the tuples of the initial integration space which determine precisely this pivot; as a rule, it is a product of homogeneous spaces.".

In the present case, the "pivot" is the 2-flat, where $dL_2$ is the density of 2-planes which is invariant under the group of rigid motions in $\mathbb{R}^d$, and the inner integration is performed in the 2-flat by using a "circumdisk representation" (pp 93-96 of [43]). Following Affentranger [2], the area elements of $L_2$ at points $A'_1, A'_2, A'_3$ are respectively $dA'_1, dA'_2, dA'_3$ and $T$ is the area of the triangle $A'_1 A'_2 A'_3$ contained in $L_2$. Polar coordinates with respect to the center $C$ of the circumcircle of these three points, of radius $\omega$, are then used to express $dA'_1 dA'_2 dA'_3$. It comes (p. 409 of [2], pp 93-96 of [43], p. 17 of [54] (1976):

$$dA'_1 dA'_2 dA'_3 = 2\omega^{-2} T dS'_1 dS'_2 dS'_3 d\omega dC \tag{14}$$

where $dS'_1, dS'_2, dS'_3$ are arc elements of the circumcircle at $A'_1, A'_2, A'_3$. Eqs 13 and 14 give then:

$$dA_1 dA_2 dA_3 = 2^{d-1} \omega^{-2} T^{d-1} dS'_1 dS'_2 dS'_3 dC dL_2 \tag{15}$$

the product $T^{d-1} dS'_1 dS'_2 dS'_3$ is then integrated over entire circumcircles because only circumcircles which are totally in the interior of the unit ball are considered [2]. The latter integral simplifies then to $\alpha(d) \omega^{2d+1}$ (eq. 6 of [2]). It remains then to express $dC$ and $dL_2$. Polar coordinates, with an origin at $O'$ $(\Delta \equiv O'C)$, are used to express $dC$: a radial factor $\delta d\delta$ appears in particular because $n = 2$ (more generally a factor $\delta^{n-1} d\delta$ for a $n$-flat) while an integration over the angular factor yields a constant. The density $dL_2$ can be represented as (eq. 9 of [2]):

$$dL_2 = h^{d-3} dh dL^*_{1[0]} dL^{(d-1)}_{2[0]} \tag{16}$$

where $h$ $(h \in (0,1))$ denotes the distance from $O$ to $L_2$ ($H \equiv OO'$, figure 1a), $dL^*_{1[0]}$ the density of oriented lines $L^*_1$ passing through $O$ and $dL^{(d-1)}_{2[0]}$ the density of 2-flats through $O'$ and perpendicular to $L^*_1$. Then, an integration over all positions of $L^*_{1[0]}$ and of $L^{(d-1)}_{2[0]}$ is performed and a factor $g(\delta, \omega, h)$ is introduced. This factor is equal to 1 if the circumcircle situated in $L_2$, with radius $\omega$ and center $C$, is contained in $B_d$ and is equal to 0 otherwise. From the previous calculations, the pdf $p_d^{(2)}(\delta, \omega, h)$ is finally deduced to be:

$$p_d^{(2)}(\delta, \omega, h) \propto g(\delta, \omega, h) \delta \omega^{2d-1} h^{d-3} \qquad (17)$$

The intersection $S_1^{(d-1)} \cap L_2$ is a circle of radius $r = \sqrt{1-h^2}$ (figure 1a). To obtain the pdf $f_d^{(2)}(\delta_r, \omega_r)$, we rewrite first eq. 17 in terms of $\Delta_r = \Delta/r$ and $\Omega_r = \Omega/r$ for which $\delta_r$ and $\omega_r$ vary both between 0 and 1:

$$p_d^{(2)}(\delta, \omega, h) d\delta d\omega dh \propto g(\delta_r + \omega_r) \delta_r \omega_r^{2d-1} d\delta_r d\omega_r \times \left( h^{d-3} \left(1 - h^2\right)^{d+1} dh \right) \qquad (18)$$

where $g(\delta_r + \omega_r)$ is equal to 1 if $\sigma_r = \delta_r + \omega_r \in (0,1)$ and 0 otherwise. The rv $H$ is independent of the pair $(\Delta_r, \Omega_r)$ and its pdf is $p_d^{(2)}(h) \propto h^{d-3}(1-h^2)^{d+1}$ $(h \in (0,1))$ in full agreement with eq. 21 [40] for $n = 2$. Eq. 18 yields finally $(T_r = 1 - \Delta_r - \Omega_r)$:

$$f_d^{(2)}(\delta_r, \omega_r) = 4d(d+1)\delta_r \omega_r^{2d-1} \quad \delta_r + \omega_r + \tau_r = 1 \quad \delta_r, \omega_r, \tau_r \in (0,1) \qquad (19)$$

The joint distribution given by eq. 19 is thus a Dirichlet distribution $Dir(2, 2d, 1)$. The latter result generalizes to:

$$f_d^{(n)}(\delta_r, \omega_r) = \frac{n(d+1)}{B(n, nd)} \times \delta_r^{n-1} \omega_r^{nd-1} \quad \delta_r + \omega_r + \tau_r = 1 \quad \delta_r, \omega_r, \tau_r \in (0,1) \qquad (20)$$

for any $d \geq 2$ and $1 \leq n \leq d$ as deduced from eqs 4.3, 4.4 and 4.6 of [3] complemented with a calculation similar to the one which leads to eq. 19. The radial part of $d\lambda_m(Q')$ which appears in eq. 4.3 of [3] $(m \equiv n)$ is $\propto \delta^{n-1} d\delta$ as mentioned above. The joint distribution given by eq. 20 is again a Dirichlet distribution $Dir(n, nd, 1)$. In addition, equations 4.3 and 4.6 of [3] show, as above, that $H$ is independent of the pair $(\Delta_r, \Omega_r)$. The distribution of $H$ is thus the same for $n$-flats selected from circumspheres of $C_d^{(n)}$ (figure 2) and for those selected from circumspheres which cut $B_d$. The latter property is required for the calculation of the bivariate pdf $p_d^{(n)}(\delta, \omega)$ from $f_d^{(n)}(\delta_r, \omega_r)$ and $w_d^{(n)}(r)$ (eq. 22).

**B. Distribution of the distance $H \equiv OO'$ and of the radius $R \equiv \sqrt{1-H^2}$ for $d \geq 2$ and $n < d$**

The pdf of $H \equiv OO'$ (figure 1a, 1d), was obtained, among others things, by Miles (theorem 4 of [44]):

$$p_d^{(n)}(h) = \frac{2}{B((d-n)/2, 1+n(d+1)/2)} \times h^{d-n-1} (1-h^2)^{\frac{n(d+1)}{2}} \qquad (21)$$

$(h \in (0,1))$ (figure 2 for $n=2$). The latter pdf was derived for $n=d-1$ by Raynaud [51]. As noted by Miles [44], $Z = H^2$ has a beta distribution, $Z \sim Be((d-n)/2, 1+n(d+1)/2)$.

The intersection of $L_n$ with the unit hypersphere $S_1^{(d-1)}$ is a hypersphere $S_R^{(n-1)}$ whose center is $O'$ and whose radius is $R \equiv O'D \equiv \sqrt{1-H^2}$ (figure 1a). It comes from eq. 21:

$$w_d^{(n)}(r) = \frac{2}{B((d-n)/2, 1+n(d+1)/2)} \times r^{n(d+1)+1} (1-r^2)^{(d-n-2)/2} \qquad (22)$$

$r \in (0,1)$, i.e. $R^2 = 1 - Z$ is beta distributed, $R^2 \sim Be(1+n(d+1)/2, (d-n)/2)$.

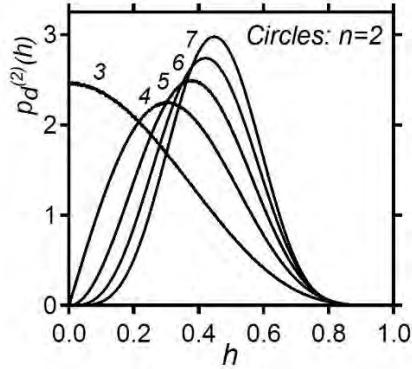

Figure 2: Circumcircles $c_d^{(2)}$ in $\mathbb{R}^d$ for $3 \leq d \leq 7$: the pdf's of the distance $H \equiv OO'$ (figures 1a and 1d) are calculated from eq. 21 and are compared with the pdf's estimated from the results of Monte-Carlo simulations in the conditions described at the end of section II.

**V. CIRCUMSPHERES $c_d^{(n)}$ OF $n+1$ POINTS CHOSEN AT RANDOM IN A UNIT $d$-BALL $(1 \leq n \leq d)$**

We derive below the joint distribution of the length $\Delta$ and of the circumradius $\Omega$, the marginal pdf's of $\Delta$, of $\Omega$ and that of their sum $\Sigma$ as well as stochastic representations of $\Delta^2$ and $\Omega^2$. The mathematical form of the joint distribution of the length $\Delta_C = OC$ and of $\Omega$ is obtained but not its explicit normalization constant. Finally, we calculate the probability $p\left(O' \notin c_d^{(n)}\right)$ that $O'$ is outside any circumsphere of $C_d^{(n)}$.

## A. Joint distribution of $\Delta$ and $\Omega$

As discussed in section IVA, the joint distribution of the triplet $(\Delta_r, \Omega_r, T_r)$ is a Dirichlet distribution $Dir(n, nd, 1)$ (eq. 20). The marginal beta distributions of $\Delta_r, \Omega_r, T_r$ and of $\Sigma_r$, are obtained from the amalgamation property of the Dirichlet distribution (end of appendix B). They are respectively $\Delta_r \sim Be(n, nd+1), \Omega_r \sim Be(nd, n+1), T_r \sim Be(1, n(d+1))$ and thus $\Sigma_r \sim Be(n(d+1), 1)$ (fig. 3 for $n = 2$). The conditional distribution of $\Delta_r$ given a circumradius $\omega_r$ is obtained directly from the Dirichlet distribution which implies that $\Delta_r/(1-\omega_r) \sim Be(n,1)$. Thus:

$$f_d^{(n)}(\delta_r | \omega_r) = \frac{n \delta_r^{n-1}}{(1-\omega_r)^n} \qquad \delta_r \in (0, 1-\omega_r) \tag{23}$$

Eq. 23 has a simple interpretation: given a circum radius $\omega_r$, the center of the circumsphere $c_d^{(n)}(\omega_r)$ is uniformly distributed in the interior of a ball of $\mathbb{R}^n$ of radius $1-\omega_r$.

To derive the joint distribution $p_d^{(n)}(\delta, \omega)$, we come back to the initial distance scale, namely $\delta = r\delta_r$ and $\omega = r\omega_r$, so that $f_d^{(n)}(\delta_r, \omega_r) \left(\delta_r + \omega_r \in (0,1)\right)$ (eq. 20) is transformed into $p_d^{(n)}(\delta, \omega | r) = (1/r^2) f_d^{(n)}(\delta/r, \omega/r)$ $(0 < \delta + \omega < r < 1)$. The latter distribution, once weighted by $w_d^{(n)}(r)$ (eq. 22), where $r$ can take any value between $\delta + \omega$ and 1 for $\delta \le \Delta \le \delta + d\delta$ and $\omega \le \Omega \le \omega + d\omega$, yields finally the desired pdf:

$$p_d^{(n)}(\delta, \omega) = \int_{\delta+\omega}^{1} p_d^{(n)}(\delta, \omega | r) w_d^{(n)}(r) dr = \int_{\delta+\omega}^{1} f_d^{(n)}\left(\frac{\delta}{r}, \frac{\omega}{r}\right) w_d^{(n)}(r) \frac{dr}{r^2}$$

$$\propto \delta^{n-1} \omega^{nd-1} \int_{\delta+\omega}^{1} 2r(1-r^2)^{\frac{d-n}{2}-1} dr \tag{24}$$

After a simple transformation of the constant factor, $p_d^{(n)}(\delta, \omega)$ is finally written as:

$$\begin{cases} p_d^{(n)}(\delta, \omega) = K_{d,n} \delta^{n-1} \omega^{nd-1} \left(1-(\delta+\omega)^2\right)^{\frac{d-n}{2}} & (\delta, \omega, \delta+\omega \in (0,1)) \\ \\ K_{d,n} = \frac{2^{n+nd}}{\sqrt{\pi}} \times \frac{\Gamma\left((n(d+1)+1)/2\right)\Gamma\left(((n+1)d+2)/2\right)}{\Gamma(n)\Gamma(nd)\Gamma\left((d-n+2)/2\right)} \end{cases} \tag{25}$$

for $1 \le n \le d$ and $d \ge 1$.

The pdf $p_d^{(n)}(\delta,\omega)$ is a polynomial in $\delta$ and in $\omega$ either when $n=d$ or when $d-n$ is even. It reduces to a constant for $n=d=1$ (eq. 4, section IIIA). If $n=d$, the joint distribution of $(\Delta,\Omega,T)$ is a Dirichlet distribution $Dir(d,d^2,1)$ (eqs 20 and 25):

$$p_d^{(d)}(\delta,\omega) = \frac{d(d+1)}{B(d,d^2)} \times \delta^{d-1}\omega^{d^2-1} \qquad (\delta,\omega,\delta+\omega \in (0,1)) \qquad (26)$$

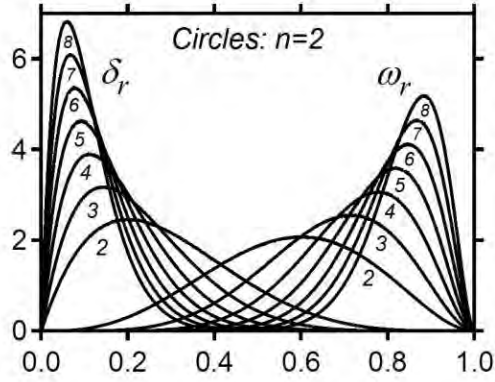

Figure 3: Circumcircles $c_d^{(2)}$ in $\mathbb{R}^d$ for $2 \leq d \leq 8$: the rescaled distance $\Delta_r$ and circumradius $\Omega_r$ are beta distributed (section IVA). Their pdf's, respectively $f_d^{(2)}(\delta_r) = 2(d+1)(2d+1) \times \delta_r(1-\delta_r)^{2d}$ and $f_d^{(2)}(\omega_r) = 2d(d+1)(2d+1) \times \omega_r^{2d-1}(1-\omega_r)^2$ $(\delta_r,\omega_r \in (0,1))$, are compared with those estimated from the results of Monte-Carlo simulations in the conditions described at the end of section II.

Any circumsphere which belongs to the $C_d^{(n)}$ family fulfils the condition $\delta+\omega \in (0,1)$. The ensuing maximum value of the sum $\delta_C + \omega = OC + CB$ is obtained from the inequality $OC \leq OO' + O'C$ in the right-angled triangle $OO'C$ and from $CB \leq CD$ (figure 1a):

$$\delta_C + \omega \leq OO' + O'D = r + \sqrt{1-r^2} \leq \sqrt{2} \qquad (27)$$

where $\delta_C, \omega, r \in (0,1)$. The maximum value of $\delta_C + \omega$ is reached for $\delta_C = \omega = r = 1/\sqrt{2}, \delta = 0$. The mathematical form of the bivariate distribution $p_d^{(n)}(\delta_C,\omega)$ is discussed in the next section.

## B. Joint distribution of $\Delta_C$ and $\Omega$

The starting pdf, valid for any $n \geq 2$ and any $d \geq n$, generalizes the one given by eq.17 for $n = 2$:

$$p_d^{(n)}(\delta,\omega,h) \propto g(\delta,\omega,h)\delta^{n-1}\omega^{nd-1}h^{d-n-1} \qquad (28)$$

Eq. 28, with $\delta,\omega,h \in (0,1)$, holds for circumspheres which belong to the family $C_d^{(n)}$, where $g(\delta,\omega,h)$ is therefore equal to 1 if $\sigma = \delta + \omega < \sqrt{1-h^2}$ and 0 otherwise. It is readily verified that eq. 28 yields, as expected, $f_d^{(n)}(\delta_r,\omega_r) \propto \delta_r^{n-1}\omega_r^{nd-1}$ (eq. 20) and $p_d^{(n)}(h) \propto h^{d-n-1}(1-h^2)^{n(d+1)/2}$ (eq. 21). As $\delta_C = \sqrt{\delta^2 - h^2}$, the pdf $p_d^{(n)}(\delta_C,\omega,h)$ writes:

$$p_d^{(n)}(\delta_C,\omega,h) \propto g(\delta_C,\omega,h)\delta_C \omega^{nd-1}(\delta_C^2 - h^2)^{(n-2)/2} h^{d-n-1} \qquad (29)$$

$\delta_C,\omega,h \in (0,1)$, where $g(\delta_C,\omega,h)$ is equal to 1 if $\omega < \sqrt{1-h^2} - \sqrt{\delta_C^2 - h^2}$ and 0 otherwise. We define $Y = H^2$ so that eq. 29 becomes:

$$p_d^{(n)}(\delta_C,\omega,y) \propto g(\delta_C,\omega,y)\delta_C \omega^{nd-1}(\delta_C^2 - y)^{(n-2)/2} y^{(d-n-2)/2} \qquad (30)$$

The equation $\omega = \sqrt{1-x} - \sqrt{\delta_C^2 - x}$, with $x \in (0,\delta_C^2)$, is first transformed into $2\omega\sqrt{\delta_C^2 - x} = 1 - \omega^2 - \delta_C^2$. It requires that $\omega^2 + \delta_C^2 \leq 1$ for a solution to exist. The latter condition implies that $\omega + \delta_C \leq \sqrt{2}$ (eq. 27). It can be obtained from a triangle like $OCB$ (figure 1a) whose angle at the vertex $C$ is $\pi - \theta$ with $\theta \in (0,\pi/2)$. It comes then, $\omega^2 + \delta_C^2 = OC^2 + CB^2 = OB^2 - 2\delta_C \omega \cos\theta \leq 1 - 2\delta_C \omega \cos\theta \leq 1$. The sought-after solution $x$ is:

$$x = \delta_C^2(1-\alpha^2), \quad \alpha^2 = (1-\omega^2-\delta_C^2)^2 / (4\omega^2 \delta_C^2) \qquad (31)$$

However, the condition $x \leq \delta_C^2$ does not guarantee that $x \geq 0$. Indeed, the numerator of $1 - \alpha^2$ is equal to $(1-\delta_C+\omega)(1-\omega+\delta_C)(1+\omega+\delta_C)(\omega+\delta_C-1)$. All factors are non-negative except the last for $\omega + \delta_C < 1$. This means that the condition $\omega < \sqrt{1-y} - \sqrt{\delta_C^2 - y}$ is fulfilled for any $y \in (0,\delta_C^2)$ for $\omega + \delta_C \leq 1$. By contrast, both $\omega + \delta_C > 1$ and $\omega^2 + \delta_C^2 \leq 1$ must be satisfied when $y \in (\delta_C^2(1-\alpha^2), \delta_C^2)$. Integrating eq. 30 with the previous conditions and defining $y = z\delta_C^2$, we obtain finally, for any $n \geq 2$ and any $d \geq n$:

$$\begin{cases} p_d^{(n)}(\delta_C,\omega) = K_d^{(n)} \delta_C^{d-1} \omega^{nd-1} \times \int\limits_{\beta(\delta_C,\omega)}^{1} z^{((d-n)/2)-1} (1-z)^{(n/2)-1} dz \\ \omega,\delta_C \in (0,1) \\ \omega^2 + \delta_C^2 \leq 1 \end{cases} : \begin{cases} = K_d^{(n)} B\left(\dfrac{d-n}{2},\dfrac{n}{2}\right) \delta_C^{d-1} \omega^{nd-1} \times \left(1 - I_{\beta(\delta_C,\omega)}\left(\dfrac{d-n}{2},\dfrac{n}{2}\right)\right) \\ \beta(\delta_C,\omega) = \left(\left(1-(\delta_C-\omega)^2\right)\left((\delta_C+\omega)^2-1\right)\bigg/\left(4\omega^2\delta_C^2\right)\right) S(\omega+\delta_C-1) \end{cases} \quad (32)$$

where $K_d^{(n)}$ is a normalization constant which has not been found explicitly in the general case, $I_{\beta(\delta_C,\omega)}\left(\dfrac{d-n}{2},\dfrac{n}{2}\right)$ is a regularized incomplete beta function (eq. 3) while $S(x)$ is a step function, i.e. $S(x)=0$ for $x \leq 0$ and $S(x)=1$ for $x>0$. The presence of $\omega$ and $\delta_C$ in the denominator of $\alpha^2$ (eq. 31) has no consequence as both variables are strictly positive as soon as $\omega+\delta_C > 1$. A Dirichlet distribution $Dir(d,d^2,1)$ is recovered from eq. 32 for $n=d$. Indeed, the sole condition which is then satisfied by circumspheres of the $C_d^{(n)}$ family is $\omega+\delta_C \leq 1$ because $\Delta \equiv \Delta_C$ (section VA). The conditional distribution $p_d^{(n)}(\delta_C|\Omega=\omega)$ varies as $\delta_C^{d-1}$ for $\delta_C \leq 1-\omega$ (eq. 32): the center $C$ is then uniformly distributed in a ball of $\mathbb{R}^d$ of radius $1-\omega$. This conclusion does not hold for $\delta_C \in \left(1-\omega, \sqrt{1-\omega^2}\right)$. Monte-Carlo simulations confirmed that conditional distributions $p_d^{(n)}(\delta_C|\Omega=\omega)$ which are estimated for $n=2$ and $d=3,..,6$ agree with those calculated with Maple from eq. 32. To conclude, a full expression of $p_d^{(n)}(\delta_C,\omega)$ has not been obtained but the latter pdf can be calculated numerically from eq. 32.

## C. Marginal distributions of $\Delta$, $\Omega$, $\Sigma$ and $\Delta_C$

For $d \geq 1$ and $1 \leq n \leq d$, the marginal distributions of $\Delta$ and of $\Omega$ are derived respectively from:

$$p_d^{(n)}(\delta) = K_{d,n} \delta^{n-1} \int\limits_0^{1-\delta} \omega^{nd-1} [1+\delta+\omega]^{(d-n)/2} [1-\delta-\omega]^{(d-n)/2} d\omega \quad (33)$$

and from:

$$p_d^{(n)}(\omega) = K_{d,n} \omega^{nd-1} \int\limits_0^{1-\omega} \delta^{n-1} [\delta+1+\omega]^{(d-n)/2} [1-\omega-\delta]^{(d-n)/2} d\delta \quad (34)$$

Both pdf's are first calculated from the following integral (eq. 3.197.8 of [28]):

$$\int_0^u x^{\nu-1}(x+\alpha)^\lambda (u-x)^{\mu-1}\,dx = \alpha^\lambda u^{\mu+\nu-1} B(\mu,\nu)\,_2F_1\left(-\lambda,\nu;\mu+\nu;-\frac{u}{\alpha}\right) \qquad \mu,\nu>0 \qquad (35)$$

then a standard transformation is used (relation 15.3.5 of [1]):

$$_2F_1(a,b;c;z) = (1-z)^{-b}\,_2F_1\left(c-a,b;c;\frac{z}{z-1}\right) \qquad (36)$$

A first representation of the marginal pdf of $\Delta$ is obtained in that way:

$$\begin{cases} p_d^{(n)}(\delta) = \dfrac{2^{(d+n)/2+nd}}{\sqrt{\pi}} \times D_{d,n} \times \delta^{n-1}(1-\delta)^{nd+(d-n)/2}\,_2F_1\left(\dfrac{d-n}{2}+1,\dfrac{n-d}{2};nd+\dfrac{d-n}{2}+1;\dfrac{1-\delta}{2}\right) \\[2mm] D_{d,n} = \dfrac{\Gamma\left((n(d+1)+1)/2\right)\Gamma\left((n+1)d/2+1\right)}{\Gamma(n)\Gamma\left(nd+(d-n)/2+1\right)} \end{cases} \qquad (37)$$

with $\delta \in (0,1)$.

A convenient quadratic transformation of the Gauss hypergeometric function $_2F_1(a,1-a;c;z)$ (eq. 15.3.32 of [1]) writes:

$$_2F_1(a,1-a;c;z) = (1-z)^{c-1}\,_2F_1\left(\frac{c-a}{2},\frac{c+a-1}{2};c;4z(1-z)\right) \qquad (38)$$

with $z<1/2$ as $z$ takes here only real values. It yields other representations of $p_d^{(n)}(\delta)$ and of $p_d^{(n)}(\omega)$ which will be used in section VI to derive stochastic representations of $\Delta$ and of $\Omega$. In the former case, it gives:

$$p_d^{(n)}(\delta) = \frac{2^n}{\sqrt{\pi}} \times D_{d,n} \times \delta^{n-1}(1-\delta^2)^{nd+(d-n)/2}\,_2F_1\left(\frac{nd}{2},\frac{nd-n+d+1}{2};nd+\frac{d-n}{2}+1;1-\delta^2\right) \qquad (39)$$

The very same method is now applied to $\Omega$ to give the two representations:

$$\begin{cases} p_d^{(n)}(\omega) = \dfrac{2^{(d+n)/2+nd}}{\sqrt{\pi}} \times W_{d,n} \times \omega^{nd-1}(1-\omega)^{(d+n)/2}\,_2F_1\left(\dfrac{d-n}{2}+1,\dfrac{n-d}{2};\dfrac{d+n}{2}+1;\dfrac{1-\omega}{2}\right) \\[2mm] W_{d,n} = \dfrac{\Gamma\left((n(d+1)+1)/2\right)\Gamma\left((n+1)d/2+1\right)}{\Gamma(nd)\Gamma\left((n+d)/2+1\right)} \end{cases} \qquad (40)$$

where $\omega \in (0,1)$ and:

$$p_d^{(n)}(\omega) = \frac{2^{nd}}{\sqrt{\pi}} \times W_{d,n} \omega^{nd-1} (1-\omega^2)^{(d+n)/2} \,_2F_1\left(\frac{n}{2}, \frac{d+1}{2}; \frac{d+n}{2}+1; 1-\omega^2\right) \tag{41}$$

The densities $p_d^{(n)}(\delta)$ and $p_d^{(n)}(\omega)$ (fig. 4 with $n=2$) are both products of powers and of a Gauss hypergeometric function. The densities given respectively by eqs 37 and 40 become polynomials when $(d-n)/2$ is either zero or a negative integer, i.e. when $n=d$ or when $d-n$ is even. Examples of pdf's $p_d^{(n)}(\omega)$ are given for circles $(n=2,\ 2 \leq d \leq 8)$ and for spheres $(n=3,\ 3 \leq d \leq 5)$ in appendix F.

The bivariate distribution of:

$$\begin{cases} \Sigma = \Delta + \Omega \\ Y = \Delta - \Omega \end{cases} \tag{42}$$

is deduced from $p_d^{(n)}(\delta, \omega)$ (eq. 25) to be:

$$\begin{cases} p_d^{(n)}(\sigma, y) = Z_{d,n} (\sigma+y)^{n-1} (\sigma-y)^{nd-1} (1-\sigma^2)^{\frac{d-n}{2}} & (\sigma \in (0,1),\ y \in (-\sigma, \sigma)) \\ Z_{d,n} = \frac{2^{n(d+1)-1}}{\sqrt{\pi}} \times \frac{\Gamma((n(d+1)+1)/2)\Gamma((n+1)d/2+1)}{\Gamma(n)\Gamma(nd)\Gamma((d-n)/2+1)} \end{cases} \tag{43}$$

Integral 3.196.3 of [28]:

$$\int_a^b (x-a)^{\mu-1}(b-x)^{\nu-1} dx = (b-a)^{\mu+\nu-1} B(\mu, \nu) \qquad \mu, \nu > 0 \tag{44}$$

is then used to calculate the marginal pdf $p_d^{(n)}(\sigma)$ by integrating $p_d^{(n)}(\sigma, y)$ over $y$ between $(-\sigma, \sigma)$:

$$p_d^{(n)}(\sigma) = \frac{2}{B(n(d+1)/2, 1+(d-n)/2)} \times \sigma^{n+nd-1}(1-\sigma^2)^{(d-n)/2} \tag{45}$$

The distribution of $\Sigma$, the sum $\Delta+\Omega$, is noteworthy as its square is simply a beta distribution, $\Sigma^2 \sim Be(n(d+1)/2, 1+(d-n)/2)$. When $n=d$, the pdf $p_d^{(d)}(\sigma)$ simplifies into $p_d^{(d)}(\sigma) = d(d+1)\sigma^{d(d+1)-1}$ in agreement with the Dirichlet distribution $Dir(d, d^2, 1)$ of $\Delta, \Omega$ and $T$ (eq. 20) which implies that $T \sim Be(1, d(d+1))$ and thus $\Sigma \sim Be(d(d+1), 1)$.

For $d \geq 2$ and $1 \leq n \leq d-1$, the square $Y_C$ of $\Delta_C$, the length $OC$, is simply expressed from Pythagoras' theorem applied to the right-angled triangle $OO'C$ (figure 1a) as:

$$Y_C \equiv \Delta_C^2 = Z + \Delta_r^2(1-Z) \tag{46}$$

where $Z = H^2$ has a beta distribution, $Z \sim Be((d-n)/2, 1+n(d+1)/2)$ (section III and Miles (1971)) and is independent of $\Delta_r \sim Be(n, nd+1)$ (section IVA). We derive now an expression of the marginal distribution $p_d^{(n)}(\delta_C)$ not directly from the bivariate distribution $p_d^{(n)}(\delta_C, \omega)$ (eq. 32) but from an application of the so-called "random variable transformation theorem" [53] to eq. 46. It writes:

$$p_d^{(n)}(y_C) \propto \int_0^1 \int_0^1 z^{(d-n)/2-1}(1-z)^{n(d+1)/2} x^{n-1}(1-x)^{nd} \delta\left[y_C - z(1-x^2) - x^2\right] dx\,dz \tag{47}$$

where $\delta\left[y_C - z(1-x^2) - x^2\right]$ is the Dirac delta function. From eq. 47, we obtain the marginal pdf of $\Delta_C$ (fig. 5, $n=2$):

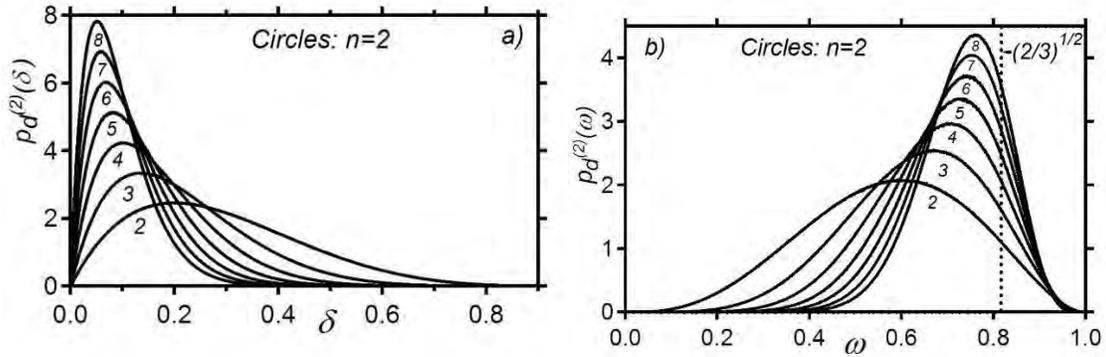

Figure 4: Circumcircles $c_d^{(2)}$ in $\mathbb{R}^d$ for $2 \leq d \leq 8$: the pdf's a) of the distance $\Delta \equiv O'C$ (figure 1a) calculated from eq. 37 and b) of the circumradius $\Omega$ calculated from eq. 40 are compared with the pdf's estimated from the results of Monte-Carlo simulations in the conditions described at the end of section II.

$$\begin{cases} p_d^{(n)}(\delta_C) = E_{d,n} \times \delta_C \left(1-\delta_C^2\right)^{n(d+1)/2} \times \int_0^{\delta_C} \frac{x^{n-1}\left(\delta_C^2 - x^2\right)^{(d-n-2)/2}(1-x)^{d(n-1)/2}}{(1+x)^{d(n+1)/2}} dx \\ E_{d,n} = \dfrac{2}{B\left((d-n)/2, n(d+1)/2+1\right) B(n, nd+1)} \end{cases} \quad (48)$$

where $\delta_C \in (0,1)$, $d \geq 2$ and $1 \leq n \leq d-1$. With the help of Maple, we verified that the pdf $p_d^{(1)}(\delta_C)$ given by eq. 48 is, as expected, identical with the pdf $p_d^{(1)}(\omega)$ given by eq. 11. Explicit expressions of $p_d^{(n)}(\delta_C)$ were not found for $n > 1$. For $n = d \geq 1$, the marginal pdf of $\Delta_C$ is the same as the marginal distribution of $\Delta_r$ (section IVA) because $O$ and $O'$ coincide. It is thus a beta distribution, $\Delta_C \sim Be(d, d^2+1)$:

$$p_d^{(d)}(\delta_C) = \frac{\delta_C^{d-1}(1-\delta_C)^{d^2}}{B(d, d^2+1)} \quad (49)$$

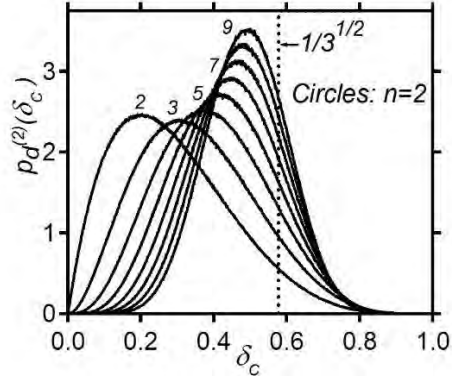

Figure 5: Circumcircles $c_d^{(2)}$ in $\mathbb{R}^d$ for $2 \leq d \leq 9$: the pdf's of the distance $\Delta_C \equiv OC$ (figure 1a) are estimated from the results of Monte-Carlo simulations in the conditions described at the end of section II.

### D. Moments of $\Delta$, $\Omega$, $\Sigma$ and $\Delta_C$

The moments of $\Delta$ and of $\Omega$ are calculated from the following integral which combines eq. 2.21.1.4 and eq. 7.3.7.8 of Prudnikov et al. (1990):

$$\int_0^1 x^{c-1}(1-x)^{b-1} {}_2F_1\left(a, 1-a; c; \frac{x}{2}\right) dx = 2^{1-b-c}\sqrt{\pi} \times \frac{\Gamma(b)\Gamma(c)}{\Gamma\left((a+b+c)/2\right)\Gamma\left((1-a+b+c)/2\right)} \quad (50)$$

We obtain:

$$E\left[\Delta^k\right]_{(d,n)} = 2^{-k} \times \frac{\Gamma(n+k)}{\Gamma(n)} \times \frac{\Gamma\left(\frac{(n+1)d+2}{2}\right)}{\Gamma\left(\frac{(n+1)d+2+k}{2}\right)} \times \frac{\Gamma\left(\frac{n(d+1)+1}{2}\right)}{\Gamma\left(\frac{n(d+1)+1+k}{2}\right)} \tag{51}$$

and

$$E\left[\Omega^k\right]_{(d,n)} = 2^{-k} \times \frac{\Gamma(nd+k)}{\Gamma(nd)} \times \frac{\Gamma\left(\frac{(n+1)d+2}{2}\right)}{\Gamma\left(\frac{(n+1)d+2+k}{2}\right)} \times \frac{\Gamma\left(\frac{n(d+1)+1}{2}\right)}{\Gamma\left(\frac{n(d+1)+1+k}{2}\right)} \tag{52}$$

The ratio of these moments is then:

$$\frac{E\left[\Omega^k\right]_{(d,n)}}{E\left[\Delta^k\right]_{(d,n)}} = \frac{(nd)_k}{(n)_k} \tag{53}$$

The duplication formula $(m)_{2p} = 4^p (m/2)_p ((m+1)/2)_p$ is used to express the even moments of $\Delta$ and of $\Omega$, $E\left[\Delta^{2p}\right]_{(d,n)}$ and $E\left[\Omega^{2p}\right]_{(d,n)}$, in terms of Pochhammer symbols:

$$E\left[\Delta^{2p}\right]_{(d,n)} = \frac{\left(\frac{n}{2}\right)_p \left(\frac{n+1}{2}\right)_p}{\left(\frac{n(d+1)+1}{2}\right)_p \left(\frac{(n+1)d+2}{2}\right)_p} \tag{54}$$

and:

$$E\left[\Omega^{2p}\right]_{(d,n)} = \frac{\left(\frac{nd}{2}\right)_p \left(\frac{nd+1}{2}\right)_p}{\left(\frac{n(d+1)+1}{2}\right)_p \left(\frac{(n+1)d+2}{2}\right)_p} \tag{55}$$

In addition, the moments of $\Sigma = \Delta + \Omega$ are obtained from eq. 45:

$$E\left[\Sigma^k\right]_{(d,n)} = \frac{\Gamma\left(\frac{(n+1)d+2}{2}\right)}{\Gamma\left(\frac{(n+1)d+2+k}{2}\right)} \times \frac{\Gamma\left(\frac{n(d+1)+k}{2}\right)}{\Gamma\left(\frac{n(d+1)}{2}\right)} \tag{56}$$

The following integral (eq. 3.211 of [28]):

$$\int_0^1 x^{\lambda-1}(1-x)^{\mu-1}(1-ux)^{-\rho}(1-vx)^{-\sigma}\,dx = B(\mu,\lambda)F_1(\lambda,\rho,\sigma;\lambda+\mu;u,v) \qquad \lambda,\mu>0 \tag{57}$$

, where $F_1(\lambda, \rho, \sigma; \lambda+\mu; u, v)$ is an Appell hypergeometric function, yields the moments of $\Delta_C$ from eq. 48 as:

$$E\left[\Delta_C^k\right]_{(d,n)} = \frac{1}{B(n, nd+1)} \times \int_0^1 x^{n-1}(1-x)^{nd} \,_2F_1\left(\frac{n(d+1)}{2}+1, -\frac{k}{2}; \frac{d(n+1)}{2}+1; 1-x^2\right) dx \tag{58}$$

Explicit expressions might be obtained for even moments from eq. 58 as $-k/2$ is then a negative integer but they become rapidly complicated as shown by the first two even moments which write:

$$\begin{cases} E\left[\Delta_C^2\right]_{(d,n)} = \dfrac{d(nd-n^2+n+1)}{(n(d+1)+1)((n+1)d+2)} \\[2mm] E\left[\Delta_C^4\right]_{(d,n)} = \dfrac{d\left(d^3n^2 - 2d^2(n^3-2n^2-2n) + d(n^4-4n^3-n^2+12n+3) - 2n^3 - 6n^2 + 8n + 6\right)}{(n(d+1)+1)(n(d+1)+3)((n+1)d+2)((n+1)d+4)} \end{cases} \tag{59}$$

The moment $E\left[\Delta_C^2\right]_{(d,n)}$ can be equally calculated from Pythagoras' theorem in the right-angled triangle $OO'C$ (figure 1a): $E\left[\Delta_C^2\right]_{(d,n)} = E\left[\Delta^2\right]_{(d,n)} + E\left[H^2\right]_{(d,n)}$ with $E\left[\Delta^2\right]_{(d,n)} = n(n+1)/\left[(n(d+1)+1)((n+1)d+2)\right]$ (eq. 54) and $E\left[H^2\right]_{(d,n)} = (d-n)/((n+1)d+2)$ (section IVB and eq. 85).

The moments $E\left[\Delta_C^k\right]_{(d,d)}$ (eq. 58) reduce to those of a beta distribution $Be(d, d^2+1)$ (eq. 49) as $_1F_0\left(-k/2;;1-x^2\right) = x^k$.

We notice finally that the moments of $1-\Delta_C^2$ can be derived from eq. 46 because $1-Z \sim Be(1+n(d+1)/2, (d-n)/2)$ and $\Delta_r \sim Be(n, nd+1)$ (section IVA) are independent. Then:

$$E\left[\left(1-\Delta_C^2\right)^k\right] = E\left[(1-Z)^k\right]E\left[\left(1-\Delta_r^2\right)^k\right] = \frac{(1+n(d+1)/2)_k}{(1+d(n+1)/2)_k} \times E\left[\left(1-\Delta_r^2\right)^k\right] \tag{60}$$

with :

$$E\left[\left(1-\Delta_r^2\right)^k\right] = \frac{1}{B(n, nd+1)} \times \int_0^1 x^{n-1}(1-x)^{nd+k}(1+x)^k \, dx \tag{61}$$

From eq. 57, we obtain finally the following explicit expression:

$$E\left[\left(1-\Delta_C^2\right)^k\right] = \frac{(1+n(d+1)/2)_k}{(1+d(n+1)/2)_k} \times \frac{(nd+1)_n}{(nd+1+k)_n} \times \,_2F_1\left(-k, n; n(d+1)+k+1; -1\right) \tag{62}$$

**E. The probability** $p\left(O' \notin c_d^{(n)}\right)$ **that** $O'$ **is not contained in** $c_d^{(n)}$

The $n$-flat $L_n$ contains almost never the center $O$ of the unit ball for any $d$ and $1 \leq n < d-1$ so that we focus on its orthogonal projection $O'$ onto $L_n$ while for $n = d$, $p\left(O' \notin c_d^{(d)}\right) \equiv p\left(O \notin c_d^{(d)}\right)$. For circumcircles, the probability $p\left(O' \notin c_d^{(2)}\right) = (1+d)/4^d$ was obtained by Affentranger [2]. The affine equivalence of $n$-flats [37], implies that it suffices to consider the rescaled pair $(\Delta_r, \Omega_r) \sim Dir(n, nd, 1)$ (eq. 20) to derive $p\left(O' \notin c_d^{(n)}\right)$ for any value of $d$ and $1 \leq n \leq d$. As $O' \notin c_d^{(n)}$ if and only if $Y_r = \Delta_r - \Omega_r > 0$, we calculate first the bivariate distribution of:

$$\begin{cases} \Sigma_r = \Delta_r + \Omega_r \\ Y_r = \Delta_r - \Omega_r \end{cases} \quad (63)$$

From the previous Dirichlet distribution, we get:

$$\begin{cases} p_d^{(n)}(\sigma_r, y_r) = J_{dn} \times (\sigma_r + y_r)^{n-1} (\sigma_r - y_r)^{nd-1} \quad \sigma_r \in (0,1), y_r \in (-\sigma_r, \sigma_r) \\ J_{dn} = \dfrac{n(d+1)}{2^{nd+n-1} B(n, nd)} \end{cases} \quad (64)$$

The marginal distribution $p_d^{(n)}(y_r)$ is then derived from the integrals $\int_{-y_r}^{1} p_d^{(n)}(\sigma_r, y_r) d\sigma_r$ for $y_r \in (-1, 0)$ and $\int_{y_r}^{1} p_d^{(n)}(\sigma_r, y_r) d\sigma_r$ for $y_r \in (0,1)$. From eq. 35, we obtain that:

$$p_d^{(n)}(y_r) = \begin{cases} \dfrac{J_{dn}}{n} \times \sum_{k=0}^{nd-1} \dfrac{(1-nd)_k}{(n+1)_k} \times (1-y_r)^{nd-1-k} (1+y_r)^k & y_r \in (-1, 0) \\ \dfrac{J_{dn}}{nd} \times \sum_{k=0}^{n-1} \dfrac{(1-n)_k}{(nd+1)_k} \times (1-y_r)^{nd+k} (1+y_r)^{n-1-k} & y_r \in (0,1) \end{cases} \quad (65)$$

Then, the probability $p\left(O' \notin c_d^{(n)}\right)$ is calculated by integrating $p_d^{(n)}(y_r)$ between 0 and 1. Expressing the sum in eq. 65 as a hypergeometric function and defining $x = (1 - y_r)/(1 + y_r)$, the latter integral becomes:

$$\dfrac{2^{nd+n} J_{dn}}{nd} \times \int_0^1 \dfrac{x^{nd}}{(1+x)^{nd+n+1}} {}_2F_1(1-n, 1; nd+1; x) dx = \dfrac{(n(d+1)-1)!}{2^{n(d+1)-1} \Gamma(n)(nd)!} \times {}_2F_1(1-n, 1; nd+1; -1) \quad (66)$$

where the right-hand side of eq. 66 is obtained first from integral 7.512.9 of [28] and then from the classical transformation $_2F_1\left(\alpha,\beta;\gamma;\dfrac{1}{2}\right)=2^\alpha\times{_2F_1}(\alpha,\gamma-\beta;\gamma;-1)$. It comes:

$$_2F_1(1-n,1;nd+1;-1)=\sum_{k=0}^{n-1}\dfrac{(-1)^k(1-n)_k}{(nd+1)_k}=\sum_{k=0}^{n-1}\dfrac{(n-k)_k}{(nd+1)_k}=\Gamma(n)\sum_{k=0}^{n-1}\dfrac{(nd)!}{(nd+k)!(n-1-k)!} \qquad (67)$$

Finally, eqs 66 and 67 yield:

$$p\left(O'\notin c_d^{(n)}\right)=\dfrac{1}{2^{n(d+1)-1}}\times\sum_{k=0}^{n-1}\binom{n(d+1)-1}{k} \qquad (68)$$

The probability $p\left(O'\notin c_d^{(n)}\right)$ (eq. 68), considered as a function of $n$, is equal to $1/2$ for $d=1$ and any $n\geq 1$. Indeed,

$$\sum_{k=n}^{2n-1}\binom{2n-1}{k}=\sum_{k'=0}^{n-1}\binom{2n-1}{n+k'}=\sum_{k'=0}^{n-1}\binom{2n-1}{n-1-k'}=\sum_{k=0}^{n-1}\binom{2n-1}{k}=2^{2n-2}.$$ The latter result cannot be interpreted in term

of the event $O'\notin c_1^{(n)}$ except for $n=1$. An unrelated explanation exists however. Indeed, the form seen in eq. 68 is found in different contexts. An example, given by Wendel [62], is fair coin tossing where the probability of at most $n-1$ heads in $N-1$ throws is given by eq. 68 in which $n(d+1)$ is replaced with $N$. The probability of at most $n-1$ heads, or that of at least $n$ tails (and vice versa), in $2n-1$ throws is thus 1/2 for any $n\geq 1$. More relevant to the present problem might be a distribution-independent result on convex hulls of $N$ i.i.d. random points in $\mathbb{R}^n$, $conv(X_1,..,X_N)$ [56]. If the distribution of these points is symmetric with respect to the origin $O''$ and assigns measure zero to every hyperplane through $O''$, then the probability that $O''\notin conv(X_1,..,X_N)$ is given by eq. 68 in which $n(d+1)$ is replaced with $N$ (eq. 1 of [62] and eq. 12.1.1 of [56]). We failed to establish an explicit relation between this problem and the present one. The former suggests nevertheless that a simplest proof of eq. 68 might exist. Some illustrative examples of $p\left(O'\notin c_d^{(n)}\right)$ are given below for $n$ ranging from 1 to 5 $(d\geq n)$:

$$\begin{cases} p\left(O'\notin c_d^{(1)}\right)=\dfrac{1}{2^d} \quad , \quad p\left(O'\notin c_d^{(2)}\right)=\dfrac{1+d}{4^d} \\ p\left(O'\notin c_d^{(3)}\right)=\dfrac{8+15d+9d^2}{8^{d+1}} \quad , \quad p\left(O'\notin c_d^{(4)}\right)=\dfrac{3+8d+9d^2+4d^3}{3\times 16^d} \\ p\left(O'\notin c_d^{(5)}\right)=\dfrac{384+1310d+2075d^2+1750d^3+625d^4}{384\times 32^d} \end{cases} \quad (d\geq n) \qquad (69)$$

The probability $p\left(O' \notin c_d^{(1)}\right)$ is derived by a simple method in section III (eq. 12). The probability $p\left(O' \notin c_d^{(2)}\right)$ agrees with the result of Affentranger [2] mentioned above. The probabilities $p\left(O' \notin c_d^{(n)}\right)$ estimated from Monte-Carlo simulations for $d \leq 5$ and $1 \leq n \leq d$ agree with those given by eq. 69. For any $n \geq 1$, $p\left(O' \notin c_d^{(n)}\right)$ is a polynomial of degree $n-1$ in $d$ divided by an exponential function of $d$ (eq. 68). It follows that $\lim\limits_{\substack{d \to \infty \\ n \text{ fixed}}} p\left(O' \notin c_d^{(n)}\right) = 0$, a result which is consistent with the fact that the center $C$ of $c_d^{(n)}$ tends to $O'$ when $d \to \infty$ (section VII).

## VI. STOCHASTIC REPRESENTATIONS OF $\Delta$ AND OF $\Omega$ FOR $d \geq 2$ AND $1 \leq n \leq d-1$

To obtain stochastic representations of $\Delta^2$ and of $\Omega^2$, we express first their pdf's. Defining $X = \Delta^2$ and $Y = \Omega^2$, they are obtained from eqs 39 and 41 to be respectively:

$$p_d^{(n)}(x) = \frac{2^{n-1}}{\sqrt{\pi}} D_{d,n} \times x^{\frac{n}{2}-1} (1-x)^{nd+(d-n)/2} {}_2F_1\left(\frac{nd}{2}, \frac{nd-n+d+1}{2}; nd + \frac{d-n}{2}+1; 1-x\right) \qquad (70)$$

with $x \in (0,1)$ and:

$$p_d^{(n)}(y) = \frac{2^{nd-1}}{\sqrt{\pi}} W_{d,n} \times y^{\frac{nd}{2}-1} (1-y)^{(d+n)/2} {}_2F_1\left(\frac{n}{2}, \frac{d+1}{2}; \frac{d+n+2}{2}; 1-y\right) \qquad (71)$$

with $y \in (0,1)$. The pdf's given by eqs 70 and 71 can be both looked at as the pdf's of products of two independent beta random variables (appendix C). Two stochastic representations both for $\Delta^2$ and for $\Omega^2$, whose parameters are calculated from eq. 92, arise from the existence of two solutions to eq. 91 as a Gauss hypergeometric function ${}_2F_1(a,b;c;x)$ remains unchanged by a permutation of $a$ and $b$. They are respectively:

$$\Delta^2 \triangleq X_1 X_2 \begin{cases} (R_1): X_1 \sim Be\left(\frac{n}{2}, \frac{nd-n+d+2}{2}\right), X_2 \sim Be\left(\frac{n+1}{2}, \frac{nd}{2}\right) \\ (R_2): X_1 \sim Be\left(\frac{n}{2}, \frac{nd+1}{2}\right), X_2 \sim Be\left(\frac{n+1}{2}, \frac{nd-n+d+1}{2}\right) \end{cases} \qquad (72)$$

and:

$$\Omega^2 \triangleq Y_1 Y_2 \begin{cases} (S_1): Y_1 \sim Be\left(\frac{nd}{2}, \frac{n+1}{2}\right), Y_2 \sim Be\left(\frac{nd+1}{2}, \frac{d+1}{2}\right) \\ (S_2): Y_1 \sim Be\left(\frac{nd}{2}, \frac{d+2}{2}\right), Y_2 \sim Be\left(\frac{nd+1}{2}, \frac{n}{2}\right) \end{cases} \qquad (73)$$

In other words, both $\Delta$ and $\Omega$ can be considered as geometric means of two independent beta random variables:

$$\begin{cases} \Delta \triangleq \sqrt{X_1 X_2} \\ \Omega \triangleq \sqrt{Y_1 Y_2} \end{cases} \tag{74}$$

We remind that $\Sigma = \Delta + \Omega \triangleq \sqrt{X}$ where $X \sim Be(n(d+1)/2, 1+(d-n)/2)$ (eq. 45).

The even moments $E[\Delta^{2p}]_{(d,n)}$ and $E[\Omega^{2p}]_{(d,n)}$, given respectively by eqs 54 and 55, are easily deduced from the corresponding stochastic representations of eqs 72 and 73 and the associated moments of beta distributions (eq. 85). We recall that $\Delta \sim Be(d, d^2+1)$, $\Omega \sim Be(d^2, d+1)$ (section IV.A) for $n = d$. These results are fully consistent with eqs 72 and 73 as shown at the end of appendix C. The stochastic representations $S_i$ $(i = 1, 2)$ (eqs 73 and 74) are convenient for Monte-Carlo simulations of separate values of $\Delta$ and $\Omega$. Unconnected values of $\Delta$ and $\Omega$ are indeed obtained in that way as the associated pair $(\Delta, \Omega)$ is not distributed as should be a pair $(\Delta, \Omega)$ whose pdf $p_d^{(n)}(\delta, \omega)$ is given by eq. 25. Finally, the stochastic representations of $\Omega$ for $n = 1$ are discussed in more detail in appendix D.

## VII. ASYMPTOTIC BEHAVIOR

We discuss briefly below the asymptotic behaviour for $d \to \infty$ while $n$ is kept fixed. When $d \to \infty$, most of the volume of a high-dimensional unit ball concentrates in a narrow annular region at its surface (see for instance Lévy (1922)). Further, all high-dimensional random vectors are almost always nearly orthogonal to each other [12]. Therefore, all pairwise distances $A_i A_j$ $(i < j, i = 1,..,j-1, j = 2,..,n+1)$ are approximately equal and all pairwise angles are approximately equal to $\pi/2$ when $d$ becomes large [29]. In other words, the lengths $OA_i$ tend to 1 and the $A_i$ 's are expected to form a standard $n$-simplex with an intervertex distance of $\sqrt{2}$ when $d \to \infty$. We notice that this holds independently of the nature of the circumsphere of the $A_i$ 's.

When $x \to \infty$, an asymptotic formula (eq. 6.1.47 of [1]), namely $\Gamma(x+a)/\Gamma(x+b) \sim x^{a-b}(1+(a-b)(a+b-1)/(2x)+O(1/x^2))$, can be applied to get the asymptotic moments of $\Delta$ and of $\Omega$ when $d \to \infty$ for a fixed value of $n$.

$$\begin{cases} E[\Delta^k]_d^{(n)} \sim O\left(\dfrac{1}{d^k}\right) \\ E[\Omega^k]_d^{(n)} \sim \left(\dfrac{n}{n+1}\right)^{k/2} + O\left(\dfrac{1}{d}\right) \end{cases} \tag{75}$$

Eq. 75 shows that the center $C$ of $c_d^{(n)}$ tends to $O'$ when $d \to \infty$. As distributions with bounded supports are uniquely determined by their moments of integer order $k$, eq. 75 indicates that the pdf of $\Omega$ transforms progressively into a narrower and narrower pdf with a peak at a position closer and closer to $(n/(n+1))^{1/2}$ when $d$ increases. This agrees with the circumradius of a standard $n$-simplex whose intervertex distance is $\sqrt{2}$ (eq. 2). The latter conclusion is consistent with the trend shown by figure 4b $(n=2)$. Figure 6b $(n=2)$ exhibits a similar trend for circumspheres of $D_d^{(n)}$ and $E_d^{(n)}$. In addition, the following limits are consistently obtained to be:

$$\begin{cases} \lim_{\substack{d \to \infty \\ n \text{ fixed}}} E\left[\Omega^k\right]_d^{(n)} = \left(\frac{n}{n+1}\right)^{k/2} \\ \lim_{\substack{d \to \infty \\ n \text{ fixed}}} E\left[\Sigma^k\right]_d^{(n)} = \left(\frac{n}{n+1}\right)^{k/2} \end{cases} \quad (76)$$

The probability $P_d^{(n)}$ $(1 \le n \le d)$ that a circumsphere belongs to the $C_d^{(n)}$ family has been derived by Affentranger [3] ($P_d^{(1)} = 1$ for any $d$). Its asymptotic value is (eq. 5.3 of [3]):

$$\lim_{\substack{d \to \infty \\ n \text{ fixed}}} P_d^{(n)} = \left(\frac{4}{n}\right)^{n/2} \times \frac{\pi^{(n-1)/2} \Gamma((n+1)/2)}{(n+1)^{(n+1)/2}} \quad (77)$$

This probability decreases rather rapidly with $n$, being less than 0.002 for $n \ge 8$. The high-dimensional properties of unit balls ensure that the asymptotic limits of $\Delta$ and of $\Omega$ for $n$ fixed and $d \to \infty$ have to be similar for the three families of circumspheres, a trend shown by Monte-Carlo simulations performed for $n = 2$.

**VIII. TAIL BEHAVIOUR OF THE CIRCUMRADIUS PDF OF ALL CIRCUMSPHERES**

The circumradius pdf $\phi_d^{(n)}(\omega)$ of all circumspheres, and not exclusively those which belong to $C_d^{(n)}$, was studied by Monte-Carlo simulations to investigate its tail behaviour. We focused more particularly on the tail behaviour of $\phi_d^{(n)}(\omega)$ for $d$ ranging from 2 to 9 and primarily for $n = d$ for a number of circumspheres $N$ most often equal to $2.10^8$. Figure 7a shows the variation with $d$ of the probabilities that a circumsphere belongs to $C_d^{(d)}$, $D_d^{(d)}$ or $E_d^{(d)}$ which are denoted respectively as $P(d)$, $Q(d)$ and $R(d)$.

The tail behaviour of $\phi_d^{(n)}(\omega)$ is actually determined by circumspheres of $E_d^{(n)}$ whose radii are unbounded while those of circumspheres of $C_d^{(n)}$ and $D_d^{(n)}$ are bounded by 1 and 2 respectively. It is more particularly related to configurations of random points $A_i (i=1,..,n+1)$ whose associated $n$-flats are about to degenerate into $(n-1)$-flats, for instance, three points $(n=2)$ which determine a plane are on the point of becoming collinear (see below). Side inequalities in any triangle $O'CM$, where $M \in L_n \cap S_1^{(d-1)} \equiv S_r^{(n-1)}$ with $O'C \equiv \delta > 1$ ($E_d^{(n)}$), $O'M = r \leq 1$, $CM = \omega$, write $\omega - 1 \leq \omega - r < \delta < \omega + r \leq \omega + 1$. Thus, the tail behavior of $\phi_d^{(n)}(\delta)$ is identical with that of $\phi_d^{(n)}(\omega)$, as observed indeed by Monte-Carlo simulations.

Monte-Carlo simulations were first performed for $n = d = 2$ to follow how moments $\langle \Omega^j \rangle, j = 1,..,12$ vary with the number $N$ of circumspheres $(N \leq 3.10^8)$. The maximum circumradius $\Omega_M(N)$ is a step function as shown by a representative simulation example in figure 7b $(N \leq 2.10^8)$ where the values of $\Omega_M(N)$ change at $N_s = 108$, 548, 903, 8527, 3834914, for instance with $\Omega_M(3834914) = 4.38\ 10^7$ for which the average deviation of points $A_1, A_2, A_3$ from a least-squares line is $\approx 10^{-8}$. The corresponding mean $\langle \Omega \rangle$ remains typically of the order of some units for any $N$ while $\langle \Omega^2 \rangle$ increases from $2.10^3$ for $N = 10^3$ to $4.10^7$ for $N = 10^8$. Most importantly, the ratios $\rho_j = \langle \Omega^{j+1}/\Omega^j \rangle$ $(j = 2,..,11)$, calculated for a given $N$, vary as $\rho_j \approx \Omega_M(N_t)$ where $N_t$ is the largest value of $N_s$ smaller or equal to $N$. They become closer and closer to $\Omega_M(N_t)$ when $j$ increases. Analogous conclusions are obtained for $3 \leq n = d \leq 9$.

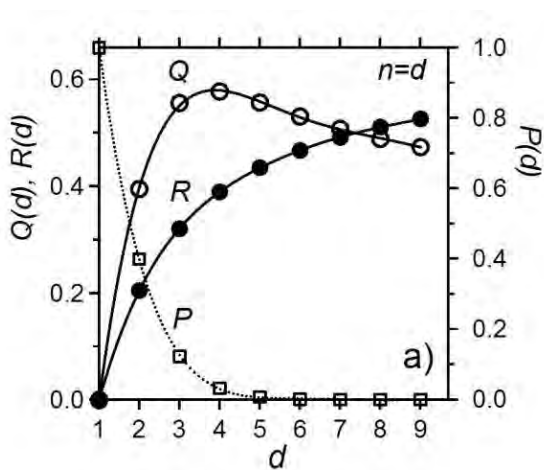
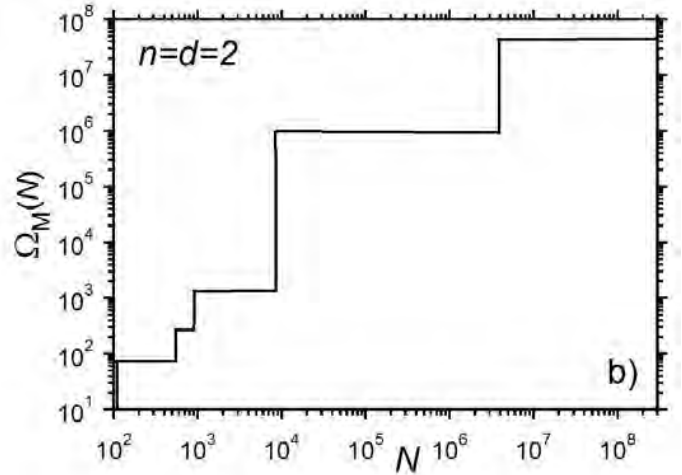

Figure 7: a) Probabilities, $P(d)$ (empty squares), $Q(d)$ (empty circles) and $R(d)$ (solid circles) that a circumsphere belongs respectively, for $n=d$, to $C_d^{(d)}$, $D_d^{(d)}$ or $E_d^{(d)}$, as estimated from Monte-Carlo simulations in the conditions described at the end of section II. Lines (solid and dotted) are guides to the eyes b) Maximum circumradius $\Omega_M(N)$ as a function of the number $N\ (\leq 2.10^8)$ of simulated circumspheres for $n=d=2$.

A working hypothesis is then that $\phi_d^{(d)}(\omega)$ has a power law tail, $\phi_d^{(d)}(\omega) \simeq a_d A_d / \omega^{a_d}$ for $\omega \to \infty$, an assumption that we considered further by studying the extreme value statistics of $\Omega$. As $\Omega$ is a positive and unbounded random variable, the relevant extreme value distribution function is the Fréchet distribution, $F(x) = \exp\left(-((x-m)/s)^{-a}\right)$, which depends on three parameters $(a,m,s)$. The associated pdf reads:

$$f(x) = \left(\frac{a}{s}\right) \times \left(\frac{x-m}{s}\right)^{-(1+a)} \exp\left(-\left(\frac{x-m}{s}\right)^{-a}\right) \qquad x > m \qquad (78)$$

It is unimodal with a mode at $m + s\left(a/(1+a)\right)^{1/a}$. It should be found for instance if the tail of the pdf of a rv $X > 0$ is $p(x) \simeq aA/x^{1+a}$ ($x \to \infty$) (see for instance [10]).

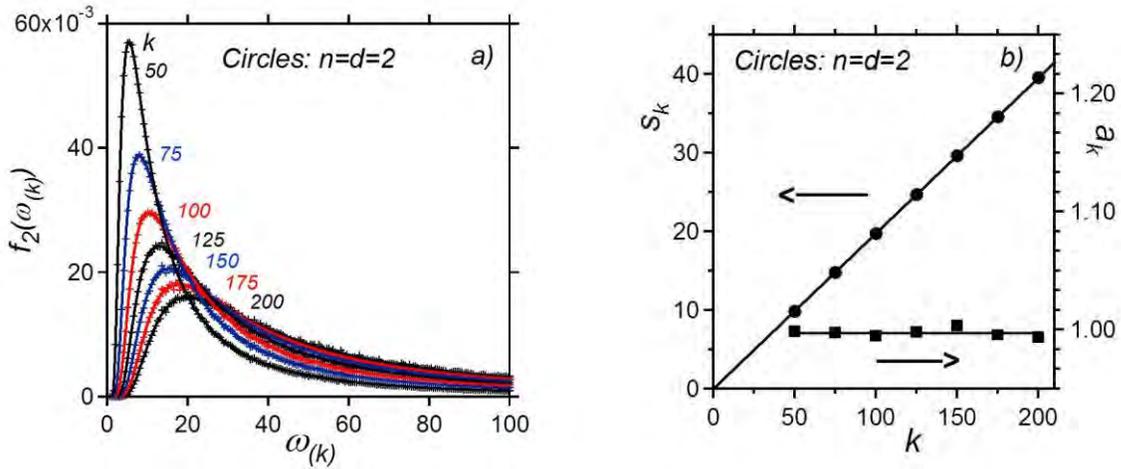

Figure 8: a) pdf's of the maximum $\omega_{(k)}$ of groups of $k$ i.i.d. circumradii of circles $(n=d=2)$ as estimated as a function of $k$ from Monte-Carlo simulations (crosses) and fitted with $f_2(\omega_{(k)})$ (eq. 79) (solid lines) b) dependence of the fitted parameters $a_k$ (right scale) and $s_k$ (left scale) on $k$.

We consider a set of $k$ i.i.d. circumradii $\Omega_j, j=1,..,k$ and we define $\Omega_{(k)} = \max(\Omega_1,..,\Omega_k)$. A limit distribution is found if there exists sequences of constants $s_k > 0$ and $m_k$ such that $\Pr\left(\left(\Omega_{(k)} - m_k\right)/s_k \leq x\right) \to F(x)$ $(x > m)$ as $k \to \infty$ ([14] p. 46). Circumradii are grouped into $N_k = N/k$ samples of size $k = 25m$ $(m=1,..,8)$ and the maximum circumradius is retained for each sample. These sets of $N_k$ maxima are then used to construct empirical pdf's of the maximum circumradius $\omega_{(k)}$. As simulated distributions are normalized in a range $(0, \omega_M)$, where $\omega_M$ is a chosen value which is not necessarily large as compared to $s_k$, they were fitted instead with the following expression:

$$f_d\left(\omega_{(k)}\right) = \left(\frac{a_k}{s_k}\right) \times \left(\frac{\omega_{(k)} - m_k}{s_k}\right)^{-(1+a_k)} \exp\left(-\left(\frac{\omega_{(k)} - m_k}{s_k}\right)^{-a_k} + \left(\frac{\omega_M - m_k}{s_k}\right)^{-a_k}\right) \quad (79)$$

$\left(\omega_{(k)} \in (m_k, \omega_M)\right)$ where a factor $\exp\left(\left((\omega_M - m)/s\right)^{-a}\right)$ has been included in eq. 79 to ensure that the calculated distribution is normalized in the same range as the simulated one. The latter factor deviates little or negligibly from one for most of the distributions of figure 8, which is shown in a range smaller than the one in which it was calculated, $(0, \omega_M = 250)$.

Figure 8 shows the pdf's of the maximum $\omega_{(k)}$ of groups of $k$ i.i.d. circumradii of circles $(n = d = 2)$ which are estimated as a function of $k$ from Monte-Carlo simulations as described above $(N = 2.10^8)$. The simulated pdf's are well fitted with $f_2\left(\omega_{(k)}\right)$ (eq. 79) for any $k$. The parameter $s_k$ varies linearly with $k$, $s_k = 0.197k$, $(0.197 \pm 0.001)$, while $a_k$ remains essentially constant, $a_k = 0.997 \pm 0.008$. The parameter $m_k$ does not show a clearly defined variation with $k$. It remains small with a mean $\langle m_k \rangle = 0.5 \pm 0.2$. If $\phi_2^{(2)}(\omega)$ is assumed to have a power-law tail with an exponent of 2, then $s_k$ should vary linearly with $k$ (eq. 1.37 of [10]). As $m_k$ is small as compared to $s_k$, the pdf mode is then expected to be equal to $s_k/2$ (here $0.099k$) (eq. 78 and below). The "experimental" mode varies indeed linearly with $k$, with a slope of $0.101 \pm 0.003$.

Results obtained for $k$ fixed and equal to 100, where $N = 2.10^8$ while $d$ $(n = d)$ varies from 2 to 9, show again that $a_{100}$ remains essentially constant, $a_{100} = 1.00 \pm 0.01$ while the parameter $s_{100}$ increases regularly with $d$ (see fig. 9 of Journal of Mathematical Physics (vol.58, no 5, 053301, 2017)).

For a given $d$, the above results lead to assume finally that $\phi_d^{(d)}(\omega) \propto 1/\omega^2$ for $\omega \to \infty$ and thus that $\phi_d^{(d)}(\delta) \propto 1/\delta^2$ for $\delta \to \infty$. Power-law tails agree also with simulation results for $d > n$, with $\phi_d^{(n)}(\omega) \propto 1/\omega^{a_d^{(n)}}$ for $\omega \to \infty$ and an exponent $a_d^{(n)}$ larger than 2 which tends to increase with $d-n$. However, we are presently unable to guess a reliable variation of it with $d$ and $n$.

## IX. CONCLUSION

We considered circumspheres which are determined almost surely by sets of $n+1$ independent random points uniformly distributed in the interior of a unit ball of center $O$ in the $d$-dimensional Euclidean space $\mathbb{R}^d$ for any $d \geq 1$ and any $1 \leq n \leq d$. The focus was put on the pair $(\Delta, \Omega)$, i.e. the distance $\Delta$ between the centre $C$ of a circumsphere and $O'$, the orthogonal projection of $O$ onto the $n$-flat determined by the $n+1$ random points, and the circumradius $\Omega$. Simple closed-form expressions of the pdf of this pair are obtained for circumspheres contained in the unit ball. Their derivation is based on previous literature results of Kingman [37], Miles [44] and Affentranger [2,3]. The $(\Delta, \Omega)$ pdf is simply a Dirichlet distribution $Dir(d, d^2, 1)$ for $n = d$. More generally, an unnoticed Dirichlet distribution, $Dir(n, nd, 1)$, plays an essential role in the calculation. Marginal pdf's of $\Delta$ and of $\Omega$, which are both products of powers and of a Gauss hypergeometric function, yield simple stochastic representations. Indeed, both random variables are distributed as geometric means of two independent beta random variables. Two pairs of different beta rv's are found both for $\Delta$ and for $\Omega$.

The mathematical form of the probability $p\left(O' \notin c_d^{(n)}\right)$ in which $n(d+1)$ is replaced with $N$, is found also in fair coin tossing where it represents the probability of at most $n-1$ heads in $N-1$ throws or it is too the probability that the origin $O''$ is outside the convex hull of $N$ i.i.d. random points in $\mathbb{R}^n$ provided that the distribution of these points is symmetric with respect to $O''$ and assigns measure zero to every hyperplane through $O''$ [56, 62].

The probability that a circumsphere cuts the unit sphere is most often larger to much larger than the probability that it does not. The problem we considered for the circumspheres of $C_d^{(n)}$ seems much simpler than it is for circumspheres of $D_d^{(n)}$ and $E_d^{(n)}$ (eq. 15 and below).

Some pdf's $\phi_d^{(n)}(\omega)$, constructed from all circumradii, irrespective of the fact whether or not they are entirely contained in the unit ball, were studied by Monte-Carlo simulations. The results obtained for $n = d$ lead to assume that the tail behaviour of $\phi_d^{(d)}(\omega)$, is a power law, $\phi_d^{(d)}(\omega) \propto 1/\omega^2$ for $\omega \to \infty$.

**ACKNOWLEDGMENTS**

I thank Dr. Richard A. Brand (Universität Duisburg-Essen) for useful discussions.

**APPENDIX A: MONTE-CARLO SIMULATIONS**

The main goal of Monte Carlo simulations is to estimate accurately all kinds of pdf's as shown in figures 2 to 9. Computations were performed with a Fortran computer code running on a standard desktop computer, valid a priori for any $d \geq 2$ and any $1 \leq n \leq d$. The rapid increase with $d$ of some exponents which appear in the expressions of pdf's (eqs 37 and 40) limits in practice the value of $d$ to $\sim 10$. Thus, simulations were carried out for $d$ ranging from 2 to 9 and a number of circumspheres $N$ of $2.10^8$. Only $\sim NP_d^{(n)}$ circumspheres are entirely contained in the unit ball $B_d$ where $P_d^{(n)}$ is given by Affentranger [3]:

$$P_d^{(n)} = \frac{(d-1)d^n \pi^{(n-1)/2}}{2(n+1)} \times \frac{B\left(\frac{n+1}{2}, \frac{nd}{2}\right) B\left(\frac{d-1}{2}, \frac{nd+1}{2}\right)}{B\left(\frac{d-n+1}{2}, \frac{nd}{2}\right)} \times \left[\frac{\Gamma\left(\frac{d}{2}\right)}{\Gamma\left(\frac{d+1}{2}\right)}\right]^{n+1} \tag{80}$$

Steps 1 and 2 described below are iterated $N$ times.

**1. Generation of sets of $n+1$ points uniformly distributed over the unit ball $B_d$ ($1 \leq n \leq d$)**

The method we used to generate a point $A$ uniformly distributed over the unit ball $B_d$ is based on the stochastic representation $A \triangleq Ru$ where $u$ is a unit vector uniformly distributed on the surface of the unit sphere $S_1^{(d-1)}$ and the pdf of $R$ is $dr^{d-1}$ ($r \in (0,1)$) (p. 75 of [22]). First, $n+1$ independent $d$-dimensional Gaussian random vectors $G_k$ ($k=1,..,n+1$), whose components are independent standard Gaussian random variables, are generated with the classical Box-Müller method [11]. Each Gaussian vector is then normalized to yield a unit vector $u_k = G_k / \|G_k\|$. The $n+1$ random vectors $u_k$ are thus independent and uniformly distributed on $S_1^{(d-1)}$ (see for instance p. 73 of [22]). Second, a number $n+1$ of independent random variables $U_k$ uniformly distributed on (0,1) are generated. Last, each unit vector $u_k$ is multiplied by $(U_k)^{1/d}$ ($k=1,..,n+1$) to obtain the sought-after set of $n+1$ points, $A_i$ ($i=1,..,n+1$), uniformly distributed over the unit ball $B_d$ of center $O$ (p. 75 of [22]). Taking arbitrarily $A_1$ as the origin of a coordinate system, we use then the Gram-Schmidt orthogonalization process to obtain an orthonormal set $(e_1,..,e_n)$ from the vectors $A_1A_k$ $(k=2,..,n+1)$.

## 2. Determination of the circumsphere characteristics ( $2 \leq n \leq d$ )

Barycentric coordinates of the center $C$ of the circumsphere of the $A_i$'s are used to write:

$$A_1C = \sum_{k=2}^{n+1} \alpha_k A_1A_k \qquad (81)$$

Then, $OC = \sum_{k=1}^{n+1} \alpha_k OA_k$ with $\alpha_1 = 1 - \sum_{k=2}^{n+1} \alpha_k$. As $C$ is equidistant from $A_1, A_2,..., A_{n+1}$, it comes:

$$A_2C.A_2C = A_3C.A_3C = .. = A_{n+1}C.A_{n+1}C = A_1C.A_1C = \omega^2 \qquad (82).$$

From these conditions, a set of $n$ linear equations in the $\alpha_k$'s is readily obtained [21]:

$$\begin{cases} A\alpha = b & (i, j = 1,..,n) \\ a_{ij} = A_1A_{i+1}.A_1A_{j+1} = a_{ji} \\ b_i = \dfrac{A_1A_{i+1}.A_1A_{i+1}}{2} \end{cases} \qquad (83)$$

Defining $S_i = OC.e_i$ $(i=1,..,n)$, the magnitudes of the two vectors $X_{//} = \sum_{k=1}^{n} S_k e_k$ and $X_\perp = OC - X_{//}$ are respectively the distance $\delta = O'C$ and the distance $h = OO'$, where $O'$ is the orthogonal projection of $O$ onto the $n$-flat obtained from the $A_i$'s (figure 1). The condition for a circumsphere to be entirely contained in the inside of $B_d$ is simply given by $\sigma = \delta + \omega < \sqrt{1-h^2}$. Identical results are found from the condition $\sigma = \delta + \omega \leq \sqrt{1-h^2}$ because circumspheres are almost never tangent to the unit sphere.

For $n=1$, the circumcenter is the middle of the segment $A_1A_2$, the radius is $\omega = \|A_1A_2\|/2$, $\delta = O'C$ (figure 1b) and $\delta_C = OC = \|OA_1 + OA_2\|/2$. For $n=d$, $\delta \equiv \delta_C = OC$ and the previous condition becomes $\sigma < 1$. Finally histograms of $\sim 1000$ bins are progressively constructed for various circumsphere characteristics.

As a quantitative test of the method described above, we compared theoretical probabilities $P_d^{(2)}$ (eq. 80, $\lim_{d\to\infty} P_d^{(2)} = \frac{\pi}{3\sqrt{3}} \approx 0.6046$) to their estimated values $P_d^{(2)}(\text{est.})$. Results are shown in table I for $d$ ranging from 2 to 9. The absolute error is of the order of some $10^{-5}$. It increases steadily when $d \geq 5$ but the agreement between simulated and calculated histograms remains very good. These conclusions hold for all studied values of $n$ and $d$.

Table I: Comparison between exact probabilities $P_d^{(2)}(\text{theor.})$ calculated from eq. 80 (Affentranger [3]) for circumscribed circles $(n=2)$ and probabilities $P_d^{(2)}(\text{est.})$ estimated from the Monte-Carlo simulations described above for $2 \leq d \leq 9$.

| $d$ | $P_d^{(2)}$ (theor.) | $P_d^{(2)}$ (est.) | $d$ | $P_d^{(2)}$ (theor.) | $P_d^{(2)}$ (est.) |
|---|---|---|---|---|---|
| 2 | $2/5 = 0.4$ | 0.40003 | 6 | $11/20 = 0.55$ | 0.55002 |
| 3 | $12\pi^2/245$ $(0.4834091..)$ | 0.48340 | 7 | $7840\pi^2/138567$ $(0.5584136..)$ | 0.55845 |
| 4 | $14/27$ $(0.5185185..)$ | 0.51851 | 8 | $494/875$ $(0.5645714..)$ | 0.56462 |
| 5 | $600\pi^2/11011$ $(0.5378042..)$ | 0.53780 | 9 | $105840\pi^2/1834963$ $(0.5692752..)$ | 0.56938 |

**APPENDIX B: THE BETA AND THE DIRICHLET DISTRIBUTIONS**

The beta distribution of a random variable $X$, $X \sim Be(q_1, q_2)$ $(q_1, q_2 > 0)$, is defined by the following probability density function [32]:

$$p_B(x) = \frac{x^{q_1-1}(1-x)^{q_2-1}}{B(q_1, q_2)} \quad x \in (0,1) \tag{84}$$

The moments of $X$ are simply given by:

$$\langle X^k \rangle = \frac{B(q_1+k, q_2)}{B(q_1, q_2)} = \frac{(q_1)_k}{(q_1+q_2)_k} \tag{85}$$

The Dirichlet distribution is a multivariate generalization of the beta distribution. It is of common use in simplices and is applied for instance in ecology [31], to model fragmentation, compositional data [5] and even in the analysis of new classifications of PNAS papers aiming to encapsulate the interdisciplinary nature of modern science [4].

The beta and the Dirichlet distributions can be simply obtained from gamma random variables. The pdf $p_G(x)$ of a gamma random variable, $X \sim \gamma(q, \theta)$, is given by:

$$p_G(x) = \frac{x^{q-1} \exp(-x/\theta)}{\theta^q \Gamma(q)} \quad (x > 0) \tag{86}$$

where $q > 0$ is the shape parameter, $\theta > 0$ the scale parameter [32]. The characteristic function of $X$ is $\varphi_X(t) = \langle e^{itX} \rangle = 1/(1 - i\theta t)^q$. A sum $S_k$ of $k$ independent gamma random variables, $\gamma(q_i, \theta)$ $(i = 1,..,k)$, with identical scale parameters and a priori different shape parameters, is itself a gamma random variable $\gamma\left(q = \sum_{i=1}^{k} q_i, \theta\right)$ as deduced from $\varphi_{S_k}(t) = \langle e^{itS_k} \rangle = 1 / \left(\prod_{j=1}^{k}(1 - it\theta)^{q_j}\right)$. As the scale parameter is irrelevant in the present context, its value is fixed at 1 from now on.

The Dirichlet distribution $Dir(\boldsymbol{q}_k)$ can be defined as follows [22]:

- First, consider a set of $k = m+1$ independent gamma random variables, $X_i \sim \gamma(q_i, 1)$ $(i = 1,...,k)$

- Second, define $S_k = \sum_{j=1}^{k} X_j$ and $L_j = X_j / S_k$ $(j = 1,..,k)$ $\left(\sum_{j=1}^{k} L_j = 1\right)$

As the sum of the $L_j$'s is equal to 1, the associated $k$-dimensional distribution is degenerate. Therefore, the Dirichlet distribution of the random vector $\boldsymbol{L}_k = (L_1, L_2,..., L_k)$ depends only on $m$ components, an arbitrary component being left aside. The pdf of the Dirichlet distribution, $\boldsymbol{L}_k \sim Dir(\boldsymbol{q}_k)$, derived from the previous definition depends on $k$ parameters collected in the vector $\boldsymbol{q}_k = (q_1, q_2,..., q_k)$. The support of the Dirichlet distribution $Dir(\boldsymbol{q}_k)$ is the unit $(k-1)$ simplex (eq. 1). The form of the pdf of the Dirichlet distribution respects the equivalence of the components of $\boldsymbol{L}_k$ (p. 17 of [22]):

$$p_D(l_1,...,l_m) = K_k \prod_{i=1}^{k} l_i^{q_i-1} \quad \left( l_k = 1 - \sum_{i=1}^{m} l_i, \quad l_i \in (0,1), \quad i=1,...,k \right) \tag{87}$$

with $K_k = \Gamma(q) / \left( \prod_{i=1}^{k} \Gamma(q_i) \right)$ and $q = \sum_{i=1}^{k} q_i$. The moment $M(r_k) = \left\langle \prod_{i=1}^{k} l_i^{r_i} \right\rangle$, $r_k = (r_1,...,r_k)$, is readily obtained to be:

$$M(r_k) = \left\langle \prod_{i=1}^{k} l_i^{r_i} \right\rangle = \frac{\prod_{i=1}^{k} (q_i)_{r_i}}{(q)_r} \tag{88}$$

where $r = \sum_{i=1}^{k} r_i$, (most of the $r_i$'s are generally taken as equal to zero). For $k=2$, the Dirichlet distribution reduces to the beta distribution (eq. 84).

The Dirichlet distribution has a notable amalgamation property [22]. If the $k$ components of $L_k \sim Dir(q_k)$ are grouped and added up into $p$ components $\Lambda_1,...,\Lambda_p$ with $\sum_{i=1}^{p} \Lambda_i = 1$, then $\Lambda_p \sim Dir(q_p^*)$ where each $q_i^*$ $(i=1,..,p)$ is the sum of the $q_j$'s corresponding to the components of the initial vector $L_k$ which add up to $\Lambda_i$. The amalgamation property is more simply understood by redefining the Dirichlet distribution $Dir(q_p^*)$ from the independent gamma rv's associated with each of the above-defined $p$ components. The latter gamma rv's, with scale parameters $q_i^*$ $(i=1,..,p)$, are themselves sums of the starting independent gamma rv's. The marginal distribution of any component $L_i$ of $L_k$ is for instance obtained by grouping the remaining $k-1$ components into a single one. Any component $L_i$ $(i=1,..,k)$ has thus a beta distribution, $L_i \sim Be(q_i, q_{\neq i})$:

$$p_{L_i}(l_i) = \frac{l_i^{q_i-1}(1-l_i)^{q_{\neq i}-1}}{B(q_i, q_{\neq i})} \quad l_i \in (0,1), \quad q_{\neq i} = \sum_{j\neq i=1}^{k} q_j \tag{89}$$

### APPENDIX C: PRODUCT OF TWO INDEPENDENT BETA RANDOM VARIABLES

A random variable $X$ is the product of two independent random variables $X_1$ and $X_2$ which are distributed according to beta distributions, respectively $X_1 \sim Be(\alpha_1, \beta_1)$ and $X_2 \sim Be(\alpha_2, \beta_2)$. Then the pdf of $X$ is given by [59]:

$$p_X(x) = \frac{B(\beta_1, \beta_2)}{B(\alpha_1, \beta_1) B(\alpha_2, \beta_2)} \times x^{\alpha_1 - 1} (1-x)^{\beta_1 + \beta_2 - 1} {}_2F_1(\alpha_1 + \beta_1 - \alpha_2, \beta_2; \beta_1 + \beta_2; 1-x) \quad (90)$$

$(x \in (0,1))$, where ${}_2F_1(a,b;c;x)$ is a Gauss hypergeometric function. By symmetry, the parameters $\alpha_1$ and $\beta_1$ may be replaced respectively by $\alpha_2$ and $\beta_2$ in eq. 90. This amounts simply to apply the Euler relation ${}_2F_1(\alpha, \beta; \gamma; x) = (1-x)^{\gamma - \alpha - \beta} {}_2F_1(\gamma - \alpha, \gamma - \beta; \gamma; x)$ [17]. A converse problem, namely the determination of the parameters of two independent beta laws from a pdf given by:

$$p_X(x) \propto x^{s-1} (1-x)^{c-1} {}_2F_1(a,b;c;1-x) \quad (91)$$

has two solutions if $a, b, c, s > 0$, $c > a$, $c > b$, $c + s > a + b$, $a \neq b$. Indeed, the identity, ${}_2F_1(a,b;c;x) \equiv {}_2F_1(b,a;c;x)$, yields two stochastic representations of $X$ as a product of two independent beta rv's, on the one hand $X_1$ and $X_2$ and on the other hand $Y_1$ and $Y_2$. These two groups of solutions are simply related by a substitution of $b$ for $a$ as seen below:

$$X \triangleq X_1 X_2 \begin{cases} X_1 \sim Be(s, c-b) \\ X_2 \sim Be(c+s-(a+b), b) \end{cases} \quad X \triangleq Y_1 Y_2 \begin{cases} Y_1 \sim Be(s, c-a) \\ Y_2 \sim Be(c+s-(a+b), a) \end{cases} \quad (92)$$

The conditions on $a, b, c, s$ given above ensure that the parameters of the beta distributions of eq. 92 are all positive. From eq. 91, the distribution of $Z = \sqrt{X}$ is:

$$p_Z(z) \propto z^{2s-1} (1-z^2)^{c-1} {}_2F_1(a,b;c;1-z^2) \quad (93)$$

For particular sets of values of $a, b, c$, for instance $(a, a+1/2, 2a+\varepsilon)$, with $\varepsilon = 0, 1$, the distribution of $Z$ becomes itself a beta distribution. To show it, we use eqs 7.3.104 $(\varepsilon = 0)$ and 7.3.105 $(\varepsilon = 1)$ of [50] to obtain:

$$ {}_2F_1(a, a+1/2; 2a+\varepsilon; 1-z^2) = \frac{2^{2a-1+\varepsilon}}{z^{1-\varepsilon}(1+z)^{2a-1+\varepsilon}} \Rightarrow Z \sim Be(2s-1+\varepsilon, 2a+\varepsilon) \quad (94)$$

Eq. 94 yields then the following explicit stochastic representations of $Z$ (eq. 92):

$$Z \triangleq \sqrt{X_1 X_2} \begin{cases} X_1 \sim Be\left(s, a+\dfrac{1}{2}\right) \\ X_2 \sim Be\left(s-\dfrac{1}{2}+\varepsilon, a+\dfrac{1}{2}\right) \end{cases} \qquad Z \triangleq \sqrt{Y_1 Y_2} \begin{cases} Y_1 \sim Be(s, a+\varepsilon) \\ Y_2 \sim Be\left(s-\dfrac{1}{2}+\varepsilon, a\right) \end{cases} \qquad (95)$$

Eq. 95 with $\varepsilon = 1$ holds for the distributions of $\Delta$ and of $\Omega$ with $n = d$ (section IV A). Eqs 70 and 93 give indeed $a = d^2/2, b = (d^2+1)/2, c = d^2+1, s = d/2$, $Z \sim Be(2s, 2a+1)$ i.e. $\Delta \sim Be(d, d^2+1)$. Similarly, from eqs 71 and 93 we get $a = d/2, b = (d+1)/2, c = d+1, s = d^2/2$, that is $\Omega \sim Be(d^2, d+1)$.

The beta rv $Z$ (eq. 95) is thus the geometric mean of two independent beta rv's. Two pairs of different beta rv's are found in general for a given $Z$. More generally, a beta rv, $Z \sim Be(\alpha, \beta)$ is distributed as the geometric mean of $n$ independent beta rv's, $X_k \sim Be((\alpha+k)/n, \beta/n), k = 0,.., n-1$: $Z \triangleq \left(\prod_{k=0}^{n-1} X_k\right)^{1/n}$ [47]. The latter property applies, for $\varepsilon = 1$, to the left part of eq. 95 with $n = 2, \alpha = 2s, \beta = 2a+1$.

## APPENDIX D. STOCHASTIC REPRESENTATIONS OF $\Omega^2$ FOR $n = 1$

For $n = 1$, the stochastic representations of the square of the half-distance between two points chosen at random in a $d$-ball ($d \geq 2$) (eq. 73) become:

$$\Omega^2 \triangleq Y_1 Y_2 \begin{cases} (S_1): Y_1 \sim Be\left(\dfrac{d}{2}, 1\right), Y_2 \sim Be\left(\dfrac{d+1}{2}, \dfrac{d+1}{2}\right) \\ (S_2): Y_1 \sim Be\left(\dfrac{d}{2}, \dfrac{d+2}{2}\right), Y_2 \sim Be\left(\dfrac{1}{2}, \dfrac{d+1}{2}\right) \end{cases} \qquad (96)$$

Lord [40] derived the distribution $p_d^{(1)}(\omega)$, as obtained from eq. 11, from the uniform distribution in the interior of the unit ball $B_d$ of the orthogonal projections on $\mathbb{R}^d$ of points uniformly distributed over the surface of a unit sphere of a $(d+2)$-dimensional Euclidean space. We profit from the pdf of the product of two independent beta rv's (eq. 90) to give a derivation of the latter property.

The spherical coordinates of a random unit vector $U$ in $\mathbb{R}^{d+2}$ depend on $d+1$ independent random angles $\Theta_1,..,\Theta_{d+1}$ where $\Theta_k \in [0, \pi), k = 1,..., d$ and $\Theta_{d+1} \in [0, 2\pi), k = 1,..., d$ (p.37 of [22], [27]). The distribution of $U$ is spherically symmetric as is the distribution of its projection $OA$ in $\mathbb{R}^d$ [40]. To characterize the latter distribution, it suffices then to consider the two first coordinates of $U$:

$$U \begin{cases} u_1 = \cos \Theta_1 \\ u_2 = \sin \Theta_1 \cos \Theta_2 \\ \ldots.. \end{cases} \tag{97}$$

The components of $\boldsymbol{OA}$ are the remaining components of $U$ in number $d$. The pdf's of $\Theta_1$ and $\Theta_2$ are respectively (p.37 of [22], [27]):

$$\begin{cases} p_1(\theta_1) = \dfrac{\sin^d \theta_1}{B(1/2,(d+1)/2)} \\ p_2(\theta_2) = \dfrac{\sin^{d-1} \theta_2}{B(1/2,d/2)} \end{cases} \tag{98}$$

From Pythagoras' theorem, the square of the length of $\boldsymbol{OA}$ is $Y = \|\boldsymbol{OA}\|^2 = 1 - (u_1^2 + u_2^2) = \sin^2 \Theta_1 \sin^2 \Theta_2$. We deduce from eq. 98 that the two independent rv's $X_1 = \sin^2 \Theta_1$ and $X_2 = \sin^2 \Theta_1$, have beta distributions, namely $X_1 \sim Be((d+1)/2, 1/2)$ and $X_2 \sim Be(d/2, 1/2)$. The pdf of their product $Y = X_1 X_2$ is then (eq. 90) $p_Y(y) \propto y^{(d-1)/2} y^{-1/2}$ as $_1F_0\left(\dfrac{1}{2}; ; 1-y\right) = y^{-1/2}$ (eq. 15.1.8 of [1]). Finally, the pdf of $R = \|\boldsymbol{OA}\| = \sqrt{Y}$ is $p_R(r) = dr^{d-1}$ $r \in (0,1)$, i.e. the point $A$ is uniformly distributed over the unit ball $B_d$.

To calculate the distribution of the half-length of the chord linking two points independently and uniformly distributed on the unit sphere $S_1^{(d+1)}$, we consider a second random unit vector $V$ in $\mathbb{R}^{d+2}$ independent of $U$. The distribution of the latter random unit vectors is invariant under any orthogonal transformation of the orthogonal group $O(d+2)$ (p. 27 of [22]). Conditioning first on $V$, we take it as the basis vector $(1,0,..,0)$ by an appropriate orthogonal transformation which depends only on $V$. The resulting distribution of $U$ remains invariant by this transformation and is independent of $V$. Therefore, the square of the sought-after half-length $Z = \Omega_{d+2}^2$ is given by:

$$Z = \frac{(U \cdot V)^2}{4} = \frac{1 - \cos \Theta_1}{2} = \sin^2\left(\frac{\Theta_1}{2}\right) \tag{99}$$

The pdf of $Z$ is obtained from the pdf $p_1(\theta_1)$ (eq. 98) which writes:

$$p_1(\theta_1) = \frac{\sin^d \theta_1}{B(1/2,(d+1)/2)} = \frac{2^d \sin^d(\theta_1/2) \cos^d(\theta_1/2)}{B(1/2,(d+1)/2)} = \frac{\sin^d(\theta_1/2) \cos^d(\theta_1/2)}{B((d+1)/2,(d+1)/2)} \tag{100}$$

$\left(2^{2q-1}/B(1/2,q)=1/B(q,q)\right)$. Then, the square $Z$ of the half-distance $\Omega_{d+2}$ has a beta distribution, $\Omega_{d+2}^2 \sim Be\left((d+1)/2,(d+1)/2\right)$, in agreement with theorem 3 of Miles [44], and the pdf of $\Omega_{d+2}$ is:

$$g_{d+2}(\varpi) = \frac{2\varpi^d \left(1-\varpi^2\right)^{(d-1)/2}}{B\left((d+1)/2,(d+1)/2\right)} \tag{101}$$

$\left(\varpi \in (0,1)\right)$, a pdf which agrees with eq. 13 of [57].

The stochastic representation $(S_1)$ of eq. 96 is thus, $\Omega^2 \triangleq Y_1 \Omega_{d+2}^2$ with $Y_1 \sim Be(d/2,1)$. If we further define $U_d^2 = Y_1$, the pdf of $U_d$ is then found to be $f(u) = du^{d-1}$ $\left(u \in (0,1)\right)$. The latter rv is simply the distance between the center of a unit ball $B_d$ and a point uniformly distributed in its interior. It may equivalently be written as $U_d \triangleq U^{1/d}$ where $U$ is uniformly distributed over $(0,1)$. Finally, representation $(S_1)$ can equally be expressed as:

$$\Omega_d^{(V)} \triangleq U^{1/d} \Omega_{d+2}^{(S)} \tag{102}$$

where the modified notation used in eq. 102 specifies the space dimensions as well as $(S)$ and $(V)$ for a uniform distribution respectively on the surface of a unit sphere and in the interior of a unit ball. As expected, the distribution of $\Omega_d^{(V)}$ is retrieved from eq. 102 to be (see eq. 11):

$$p_d^{(1)}(\omega) = \int_\omega^1 f\left(\frac{\omega}{x}\right) g_{d+2}(x) \frac{dx}{x} \propto \omega^{d-1} \int_\omega^1 \left(1-x^2\right)^{(d-1)/2} dx \propto \omega^{d-1} \int_0^{1-\omega^2} y^{(d-1)/2} (1-y)^{-1/2} dy$$

$$\propto \omega^{d-1} I_{1-\omega^2}\left((d+1)/2, 1/2\right) \tag{103}$$

**APPENDIX E:** see *Journal of Mathematical Physics (vol.58, no 5, 053301, (2017))* or *hal-01517928 (https://hal.archives-ouvertes.fr/)*

**APPENDIX F: EXAMPLES OF CIRCUMRADIUS PDF's OF CIRCLES $(n=2)$, OF SPHERES $(n=3)$ AND OF HYPERSPHERES $(n=d)$ $(\omega \in (0,1))$**

For circles, $n=2$, with $2 \leq d \leq 8$, the pdf's $p_d^{(2)}(\omega)$ are:

$$p_2^{(2)}(\omega) = 60\omega^3 (1-\omega)^2$$

$$p_3^{(2)}(\omega) = \frac{2205}{16} \times \omega^5 \times \left[\left(2+\omega^2\right)\sqrt{1-\omega^2} - 3\omega \times \arccos(\omega)\right]$$

$$p_4^{(2)}(\omega) = 360\omega^7 (1-\omega)^3 (3+\omega)$$

$$p_5^{(2)}(\omega) = \frac{495495}{1024} \times \omega^9 \times \left[ \left(8 + 9\omega^2 - 2\omega^4\right)\sqrt{1-\omega^2} - 15\omega \times \arccos(\omega) \right]$$

$$p_6^{(2)}(\omega) = \frac{13104}{5} \times \omega^{11} (1-\omega)^4 \left(5 + 4\omega + \omega^2\right)$$

$$p_7^{(2)}(\omega) = \frac{14549535}{16384} \times \omega^{13} \left[ \left(48 + 87\omega^2 - 38\omega^4 + 8\omega^6\right)\sqrt{1-\omega^2} - 105\omega \times \arccos(\omega) \right]$$

$$p_8^{(2)}(\omega) = \frac{26928}{7} \times \omega^{15} (1-\omega)^5 \left(35 + 47\omega + 25\omega^2 + 5\omega^3\right)$$

For spheres, $n = 3$, with $3 \leq d \leq 5$, the pdf's $p_d^{(3)}(\omega)$ write:

$$p_3^{(3)}(\omega) = 1980\omega^8 (1-\omega)^3$$

$$p_4^{(3)}(\omega) = \frac{229376}{11\pi} \times \omega^{11} \times \left[ \left(1 + 4\omega^2\right)\arcsin(\omega) - \frac{\omega}{3}\sqrt{1-\omega^2}\left(13 + 2\omega^2\right) \right]$$

$$p_5^{(3)}(\omega) = 12240\omega^{14} (1-\omega)^4 (4+\omega)$$

For hyperspheres, $n = d$ $(d \geq 1)$ (eq. 35):

$$p_d^{(d)}(\omega) = \frac{\omega^{d^2-1}(1-\omega)^d}{B(d^2, d+1)}$$